\documentclass{amsart}
\usepackage{amssymb,mathrsfs,amsmath,dsfont,array}% Using array, you can center the heading line of a tabular.
\usepackage{multirow,arydshln}
\usepackage{tikz}
\usetikzlibrary{arrows,positioning}
\usepackage{capt-of}% It helps you to use caption for tikzpictures as well.
\usepackage{mathbbol}     % Using these you can use mathbb fonts for lowercase characters as well.
\usepackage{color}
\usepackage{verbatim} %Allows use of \begin{comment} ... \end{comment}
\usepackage[all,cmtip]{xy}
\usepackage{amsbsy}% this helps us to use bold fonts for mathematical symbols as well. We then use $\pmb{any symbol}.$
\usepackage{amsmath}% this allows \tag{any symbol} to be used  for the equations.
\usepackage{nccmath}% this allows you to have equation with position on the left. \begin{fleqn}\begin{flalign}\end{fleqn}\end{fleqn}

\usepackage{hyperref}

\newtheorem{Thm}{Theorem}[section]
\newtheorem{Pro}[Thm]{Proposition}
\newtheorem{cor}[Thm]{Corollary}
\newtheorem{lem}[Thm]{Lemma}

\newtheorem{deft}[Thm]{Definition}

\numberwithin{equation}{section}
\renewcommand{\tilde}{\widetilde}

\def\a1s{a_1,\cdots, a_s}
\def\a{\alpha}

\def\andd{\quad\hbox{and}\quad}

\def\ad{\hbox{ad}}

\def\b{\beta}

\def\bl4{B_{\ell\geq4}}

\def\bbbc{{\mathbb C}}

\def\d{\delta}
\def\D{\Delta}

\def\gG{\mathscr{G}}
\def\gg{{\mathcal G}}

\def\fg{\mathfrak{g}}

\def\hH{\mathscr{H}}

\def\hh{{\mathcal H}}

\def\fh{\mathfrak{h}}

\def\ep{\epsilon}
\def\fm{(\cdot,\cdot)}

\def\bbbr{{\mathbb R}}

\def\1k{\frac{1}{k}}
\def\op{\oplus}
\def\ot{\otimes}

\def\sub{\subseteq}
\def\sg{\sigma}

\def\pf{\noindent{\bf Proof. }}

\def\supp{\hbox{\rm supp}}
\def\sspan{\hbox{\rm span}}

\def\bbbz{{\mathbb Z}}

\def\1il{1\leq i\leq\ell}

\def\sgn{\hbox{sgn}}

\begin{document}
\title{Classification of bases of twisted affine root supersystems}
\maketitle
\centerline{ Malihe
Yousofzadeh\footnote{Department of Mathematics, University of Isfahan, Isfahan, Iran,
P.O.Box 81745-163, and School of Mathematics, Institute for Research in
Fundamental
Sciences (IPM), P.O. Box: 19395-5746, Tehran, Iran.
Email address: ma.yousofzadeh@sci.ui.ac.ir \& ma.yousofzadeh@ipm.ir.\\
This research  was in part supported by
a grant from IPM (No. 98170424) and  is partially carried out in
IPM-Isfahan Branch.}}
\vspace{2cm}
\begin{abstract}Following the definition of a root basis of an affine root system, we define  a base of the root system of an affine  Lie superalgebra to be a  linearly independent subset $B$ of its root system  such that  each root  can be written as a linear combination of elements of $B$ with integral  coefficients such that all coefficients are nonnegative or all coefficients are nonpositive.
Characterization and classification of  bases of  root systems of affine Lie algebras are known in the literature; in fact,  up to $\pm 1$-multiple, each base of an affine root system is conjugate with the standard base under  the Weyl group action.
In the super case, the existence of those self-orthogonal roots which are not orthogonal  to at least one other root, makes the situation more complicated. In this work, we give a complete characterization of bases of a twisted affine root supersystem. We precisely describe  and classify them.
%\keywords{Twisted Affine Lie Superalgebra\and  Affine root supersystem\and  Base\and Quasi-Weyl group.}
% \PACS{PACS code1 \and PACS code2 \and more}
%\subclass{17B67}
\end{abstract}

\section{Introduction}
The concept of a base is crucial to study highest weight modules of finite dimensional basic classical simple Lie superalgebras and affine Lie (super)algebras.
The purpose of this work is finding a  characterization and giving  the classification  of bases of the root systems of twisted affine Lie superalgebras; affine Lie superalgebras are the super version of affine Lie algebras which have been  defined, constructed and classified  by J.W. Van de Lour \cite{van-thes} in 1986. In what follows we briefly recall  their structure and classification:

Suppose that $n$ is a positive integer and $I=\{1,\ldots,n\}.$
Let $\tau$ be a subset of $I$ and $A$ be a nonzero indecomposable  symmetrizable  $n\times n$-matrix with complex  entries (in which by symmetrizable, we mean that  $A$ has   a decomposition $A=DB$ with an invertible diagonal matrix  $D$ and  a symmetric matrix $B$), satisfying the following:
$$\begin{array}{ll}
\hbox{if $a_{i,j}=0,$}&\hbox{then $a_{j,i}=0;$}\\
\hbox{if $a_{i,i}=0,$}&\hbox{then $i\in\tau;$}\\
\hbox{if $a_{i,i}\neq0,$}&\hbox{then $a_{i,i}=2;$}\\
\hbox{if $a_{i,i}\neq 0,$} &\hbox{then  $a_{i,j}$ (resp. $a_{i,j}/2$) is a nonpositive integer}\\
&\hbox{for $i\in I\setminus \tau$ (resp. $i\in \tau$) with $i\neq j.$}
\end{array}$$
Fix a complex vector space $\hh$ of dimension $n+corank(A)$ and denote its dual space with $\hh^*.$ Then  there exist linearly independent subsets $$\Pi:=\{\a_i\mid i\in I\}\sub \hh^* \andd \check\Pi:=\{\check\a_i\mid i\in I\}\sub \hh$$ such that $$\a_j(\check\a_i)=a_{i,j}\quad\quad (i,j\in I).$$ Let  $\tilde\gg(A,\tau)$ be the Lie superalgebra generated by $\{e_i,f_i\mid i\in I\}\cup \hh$ subject to the following relations:
\begin{equation*}
\begin{array}{ll}
\;[e_i,f_j]=\d_{i,j}\check\a_i, & [h,h']=0,\\
\;[h,e_j]=\a_j(h)e_j,& [h,f_j]=-\a_j(h)f_j,\\
{\rm deg}(h)=0, & {\rm deg}(e_i)={\rm deg}(f_i)=\left\{\begin{array}{ll}
0&\hbox{if  $i\not\in \tau$}\\
1&\hbox{if $i\in\tau$}
\end{array}\right.
\end{array}
\end{equation*}
for $i,j\in I$ and $h,h'\in\hh,$ where $\d_{i,j}$ denotes the Kronecker delta. There is a unique maximal ideal $\mathfrak{i}$ of $\tilde\gg(A,\tau)$ intersecting $\hh$ trivially. The Lie superalgebra $\gg:=\tilde\gg(A,\tau)/\frak{i}$  is called  an {\it affine Lie superalgebra} if it is not of finite dimension but  of finite {\it growth}; see \cite[\S 6.1]{van-thes}.
 By a convention, the images of $e_i,f_i$ and $h$ $(i\in I, h\in\hh)$ in $\gg(A,\tau)$ under the canonical projection map are  denoted again by $e_i,f_i$ and $h,$ respectively. Then $\gg$ is generated by $\{e_i,f_i\mid 1\leq i\leq n\}\cup \hh.$
The Lie superalgebra $\gg$ has a root space decomposition with respect to $\hh.$ Moreover, $\Pi$ is a subset of the corresponding root system $R$ and each element  of $R$ is written as a $\bbbz$-linear combination  of elements of $\Pi$ whose coefficients are all nonnegative or all are nonpositive; we refer to $\Pi$ as the standard base of $R.$ Elements of $R$ are called roots and roots corresponding to even and odd part of $\gg$ are called respectively even and odd.

  As for affine Lie algebras, affine Lie superalgebras are constructed using an affinization process:
Suppose that $\fg:=\fg_0\op\fg_1$ is a finite dimensional  basic classical simple Lie superalgebra with a Cartan subalgebra $\fh\sub\fg_0.$ Suppose that $\kappa$ is  a nondegenerate supersymmetric invariant  even  bilinear form on $\fg$ and $\sg$ is an automorphism of $\fg$ of order $k.$ Since $\sg$ preserves $\fg_0$ as well as $\fg_1,$ we have
\begin{align*}
\fg_i=\bigoplus_{s=0}^{k-1}{}^{[s]}\fg_i\quad\hbox{where}\quad {}^{[s]}\fg_i=\{x\in\fg_i\mid \sg(x)=\zeta^sx\}\quad\quad(i\in\bbbz_2,\; 0\leq s\leq k-1)
\end{align*}
in which $\zeta$ is the $k$-th primitive root of unity.  Then
\begin{equation*}\label{aff}
\widehat\fg:=\widehat\fg_0\op\widehat\fg_1\quad\hbox{where}\quad\widehat\fg_i:=\bigoplus_{s=0}^{k-1}({}^{[s]}\fg_i\ot t^s\bbbc[t^{\pm  k}])\quad\quad(i\in\bbbz_2)
\end{equation*}
is a subalgebra of the current superalgebra $\fg\ot\bbbc[t^{\pm1}].$ Set $$\gG=\bigoplus_{s=0}^{k-1}({}^{[s]}\fg\ot t^s\bbbc[t^{\pm k}])\op\bbbc c\op\bbbc d\andd \mathscr{H}:=(({}^{[0]}\fg\cap\fh)\ot 1)\op\bbbc c\op\bbbc d,$$ then $\gG$ together with
$$[x\ot t^p+rc+sd,y\ot t^q+r'c+s'd]:=[x,y]\ot t^{p+q}+p\kappa(x,y)\d_{p+q,0}c+sqy\ot t^q-s'px\ot t^p$$ is  an {\it affine Lie superalgebra} and $\mathscr{H}$ is a Cartan subalgebra of $\gG.$  The Lie superalgebra $\gG$ is denoted by $X^{(k)}$ where $X$ is the type of $\fg$ and $X\neq A(\ell,\ell)$. The definition of $A(\ell,\ell)^{(1)}$ is a little bit different. The Lie superalgebra $X^{(k)}$  is called {\it twisted} if $k>1$ and it is called {\it non-twisted} if $k=1.$ Twisted and non-twisted affine Lie superalgebras cover all affine Lie superalgebras. We refer to the root system $R$ of $\gG,$ an affine root supersystem.

 Define $\d$ to be the functional on $\hH$ vanishing  on $(({}^{[0]}\fg\cap\fh)\ot 1)\op\bbbc c$ and  mapping $d$ to $1.$
If $\dot R$ is the  root system of $\fg,$ $\dot R+\bbbz\d$ is the root system of   $\gG=X^{(1)}.$

Using the form $\kappa$ on $\fg,$ one has a   nondegenerate supersymmetric invariant even bilinear form on $\gG$ which in turn, induces a symmetric bilinear form $\fm$ on the dual space $\hH^*$ of $\hH.$ Moreover, if $\gG$ is twisted, for nonnegative integers  $m$ and $n$ as in the following table and $V:=\sspan_\bbbr R$, there are  $\ep_i,\d_p\in V$ $(1\leq i\leq m,1\leq p\leq n)$ such that
$$(\ep_i,\ep_j):=\d_{i,j},\;\; (\d_p,\d_q):=-\d_{p,q},\;\; (\ep_i,\d_p):=0,\;\;(\d,V)=\{0\},$$ $\gG$ is one of   the Lie superalgebras in the first column of the following table  and $R$ is correspondingly  expressed as follows:

\begin{table}[h]\caption{} \label{table1}
% "h" means that the table will be placed "here", if we use [t], it will be placed on "top" and [b] indicates the "below". If we do not use this command, the table is floating.
 {\footnotesize \begin{tabular}{|c|l|}
\hline
$\hbox{Type}$ &\hspace{3.25cm}$R$ \\
\hline
$\stackrel{A(2m,2n-1)^{(2)}}{\hbox{\tiny$(n\neq 0)$}}$&$\begin{array}{rcl}
\bbbz\d
&\cup& \bbbz\d\pm\{\ep_i,\d_j,\ep_i\pm\ep_r,\d_j\pm\d_s,\ep_i\pm\d_j\mid i\neq r,j\neq s\}\\
&\cup& (2\bbbz+1)\d\pm\{2\ep_i\mid 1\leq i\leq m\}\\
&\cup& 2\bbbz\d\pm\{2\d_j\mid 1\leq j\leq n\}
\end{array}$\\
\hline
$\stackrel{A(2m-1,2n-1)^{(2)}}{\hbox{\tiny$m,n\neq0,(m,n)\neq (1,1)$}}$& $\begin{array}{rcl}
\bbbz\d&\cup& \bbbz\d\pm\{\ep_i\pm\ep_r,\d_j\pm\d_s,\d_j\pm\ep_i\mid i\neq r,j\neq s\}\\
&\cup& (2\bbbz+1)\d\pm\{2\ep_i\mid 1\leq i\leq m\}\\
&\cup& 2\bbbz\d\pm\{2\d_j\mid 1\leq j\leq n\}
\end{array}$\\
\hline
$\stackrel{A(2m,2n)^{(4)}}{\hbox{\tiny$(n,m)\neq (0,0)$}}$& $\begin{array}{rcl}
\bbbz\d&\cup&  \bbbz\d\pm\{\ep_i,\d_j\mid 1\leq i\leq m,\;1\leq j\leq n\}\\
&\cup& 2\bbbz\d\pm\{\ep_i\pm\ep_r,\d_j\pm\d_s,\d_j\pm\ep_i\mid i\neq r,j\neq s\}\\
&\cup&(4\bbbz+2)\d\pm\{2\ep_i\mid 1\leq i\leq m\}\\
&\cup& 4\bbbz\d\pm\{2\d_j\mid 1\leq j\leq n\}
\end{array}$\\
\hline
$\stackrel{D(m+1,n)^{(2)}}{\hbox{\tiny$(n\neq 0)$}}$& $\begin{array}{rcl}
\bbbz\d&\cup&  \bbbz\d\pm\{\ep_i,\d_j\mid 1\leq i\leq m,\;1\leq j\leq n\}\\
&\cup& 2\bbbz\d\pm\{2\d_j,\ep_i\pm\ep_r,\d_j\pm\d_s,\d_j\pm\ep_i\mid i\neq r,j\neq s\}
\end{array}$\\
\hline
 \end{tabular}}
 \end{table}\noindent where  we make a  convention that if  $m=0$ or $n=0,$ then their  corresponding expressions disappear. We set
 $$
\begin{array}{ll}
R^0 := \bbbz\d ,&\quad
R ^\times := R \setminus R ^0 ,\\
R ^\times_{re} := \{\alpha \in R \mid (\alpha,\alpha) \neq 0\},&\quad
R _{re} := R ^\times_{re} \cup \{0\}\quad(\hbox{real roots}),\\
R ^\times_{ns} :=\{\alpha \in R^\times | (\alpha,\alpha) = 0\},&\quad
R_{ ns} := R^ \times_{ns} \cup \{0\}\quad(\hbox{nonsingular roots})
\end{array}
$$
and mention that the existence of  nonzero nonsingular roots is a phenomena which does not occur for  affine Lie algebras.

We recall that the Weyl group of an affine root (super)system $R$ is the group generated by  linear transformations $\sg_\a,$ for  nonzero real roots $\a,$ on $V=\sspan_\bbbr R$ mapping $\b\in V$ to $\b-2\frac{(\b,\a)}{(\a,\a)}\a.$ As all nonzero nonsingular roots are odd and two times of an odd real root is an even real root  we call $\sg_\a$'s ($\a\in R_{re}^\times$) even reflections.

Following the definition of a root basis of an affine root system; see \cite[\S 5.9]{K2}, we define  a base of the root system $R$ of an affine  Lie superalgebra $\gG$ to be a  linearly independent subset $B$ of $R$ such that  each element of $R$  can be written as a linear combination of elements of $B$ with integral  coefficients such that all coefficients are nonnegative or all coefficients are nonpositive.

Characterization and classification of  bases of  root systems of affine Lie algebras are known in the literature; in fact,  up to $\pm 1$-multiple, each base of an affine root system is conjugate with the standard base under  the Weyl group action; see \cite[Pro. 5.9]{K2}.
The existence of nonsingular roots in the super case makes the situation more complicated.

In \cite{serg2},  V.~Serganova defines a base\footnote{She defines a base for all Kac-Moody Lie superalgebras in general.} of $R$ to be a linearly independent subset $B$ of $R$ such that for each $\a\in B,$ there are root vectors $x_\a$ and $y_\a$ corresponding to $\a$ and $-\a$ respectively such that $\{x_\a,y_\a\mid \a\in B\}\cup \hH$ generates the affine Lie superalgebra $\gG$ and $[x_\a,y_\b]=0$ for $\a\neq \b$; we refer to such a base as an S-base. It is trivial that each S-base is a base in our sense.

One of the most important differences between bases and S-bases is that bases are purely combinatorial  objects  while the algebraic feature  of an affine Lie superalgebra gets  involved to  define  S-bases.

The author in \cite{serg2}  introduces ``odd reflections'' to get a similar result as in non-super case to describe S-bases; more precisely, for a nonsingular root $\a$ of an S-base $ B ,$ she defines a map
\[s_\a: B \longrightarrow R\quad\quad \b\mapsto\left\{
\begin{array}{cl}
-\a& \hbox{if $\b=\a$}\\
\a+\b& \hbox{if $\b+\a\in R$}\\
\b & \hbox{otherwise}
\end{array}
\right.
\] called an odd reflection and shows that for an S-base $ B $ and a nonsingular root $\a\in B ,$ $s_\a( B )$ is also an S-base of $R,$ moreover, up to $\pm1$-multiple, each S-base of $R$ is obtained  from $ B $ by iterating the  action of odd and even reflections  \cite[Thm.~8.3 \& Lem.~2.4]{serg2}.

In contrast to an even reflection which depends only on a nonzero real root, an odd reflection based on a nonzero nonsingular root,  depends also on an S-base to which  the root  belongs., i.e.,  to get all S-bases, one needs to choose an S-base $\Pi$ and form all S-bases obtained from $\Pi$ by even reflections as well as  odd reflections based on nonsingular roots of $\Pi$ and then form  all obtained S-bases from each of these new  S-bases and continue this process.

In this work, we give a complete characterization of bases of a twisted affine root supersystem $R.$ We precisely describe  and classify them. As a by-product, we show that for twisted affine root supersystems, there is no difference between bases and S-bases.

To have a better view of what we have in this paper, let us state our main theorem: Suppose that $R$ is a twisted affine root supersystem (see Table \ref{table1}) and keep the same notation as above. We
define
\begin{align}\tag{$\dagger$}\label{form*}
\fm_{*}:&V\times V\longrightarrow \bbbc\\
&(\d,V)_*:=\{0\},\; (\ep_i,\ep_j)_*:=\d_{i,j},\; (\d_p,\d_q)_*:=\d_{p,q}\andd (\ep_i,\d_p)_*=0\nonumber
\end{align} for $1\leq i,j\leq m$ and $1\leq p,q\leq n.$
For $\a\in R$ with $(\a,\a)_*\neq 0,$ the {\it quasi-reflection} $r_\a$ (based on $\a$) is defined to be the linear  automorphism of $V$ mapping $a\in V$ to
$a-2\frac{(a,\a)_*}{(\a,\a)_*}\a.$ The subgroup  $W$ of the group of automorphisms of $V$ generated by
$\{r_\a\mid \a\in R;\; (\a,\a)_*\neq 0\}$ is called the {\it quasi-Weyl group}
of $R.$

\medskip

For $v=\sum_{i=1}^mt_{\ep_i}\ep_i+\sum_{p=1}^nt_{\d_p}\d_p+m\d\in V,$ we define the {\it support} of $v$ to be
\[\supp(v):=\{\ep_i\mid 1\leq i\leq m,\; t_{\ep_i}\neq 0\}\cup \{\d_p\mid 1\leq p\leq n,\; t_{\d_p}\neq 0\}.\]
For $\zeta\in\supp(v),$ we set \[\sgn(\zeta;v)=\left\{\begin{array}{ll}
1& \hbox{if $t_\zeta>0,$}\\
-1& \hbox{if $t_\zeta<0$}
\end{array}
\right.
\]
and for $S\sub V,$ we set
\[\supp(S):=\bigcup_{\a\in S}\supp(\a).\]

\begin{Thm}[Main Theorem] \label{main}
Set $\ell:=m+n.$ Assume  $k_1,\ldots,k_\ell\in\bbbz$  and  $\zeta_1,\ldots\zeta_\ell\in\{\pm \ep_i,\pm\d_p\mid 1\leq i\leq m,1\leq p\leq n\}$ with $\zeta_i\neq \pm\zeta_j$ for $i\neq j.$ Set   $$\theta_i:=\zeta_i+k_i\d\quad\quad (1\leq i\leq \ell)$$ and assume $\Pi$ is as in the following table:

{\small
\begin{table}[h]\caption{}\label{table2}
\centering
\begin{tabular}{|c|>{\raggedright\arraybackslash}m{59mm}|c|}
\hline
$R$&\multicolumn{1}{>{\centering\arraybackslash}m{59mm}|}{$\Pi$}& \hbox{Further Assumptions}
\\\hline
{\footnotesize$A(2m,2n)^{(4)}$}&$\{\theta_i-\theta_{i+1},\theta_\ell,\d-\theta_1\}_{i=1}^{\ell-1}$&$k_i\equiv k_j\;\hbox{(mod 2)}$\\\hline
{\footnotesize $D(m+1,n)^{(2)}$}&$\{\theta_i-\theta_{i+1},\theta_\ell,\d-\theta_1\}_{i=1}^{\ell-1}$&$k_i\equiv k_j\;\hbox{(mod 2)}$\\
\hline
\multirow{4}{*}{\footnotesize{$ A(2m,2n-1)^{(2)}$}} &\multirow{2}{*}{$\{-2\theta_1+\d,\theta_i-\theta_{i+1},\theta_\ell\}_{i=1}^{\ell-1}$}& $\supp(\theta_1)\sub\{\ep_i\mid 1\leq i\leq m\}$ \\
 &&\\\cdashline{2-3}
 &\multirow{2}{*}{$\{\d-(\theta_1+\theta_2),\theta_i-\theta_{i+1},\theta_\ell\}_{i=1}^{ \ell-1}$}& $\supp(\theta_1)\sub\{\d_p\mid 1\leq p\leq n\}$ \\
 &&\\
\hline
\multirow{16}{*}{\footnotesize$\begin{array}{c}A(2m-1,2n-1)^{(2)}\\ \hbox{\footnotesize$(m,n)\neq (1,1)$}\end{array}$}
 &\multirow{4}{*}{$\{\d-(\theta_1+\theta_2),\theta_i-\theta_{i+1},\theta_{\ell-1}+\theta_\ell\}_{i=1}^{\ell-1}$}& $\supp(\theta_1)\sub\{\d_p\mid 1\leq p\leq n\}$ \\
&&\\
&&$\supp(\theta_\ell)\sub\{\ep_i\mid 1\leq i\leq m\}$\\
 &&\\
\cdashline{2-3}
 & \multirow{4}{*}{$\{-2\theta_1,\theta_i-\theta_{i+1},\theta_{\ell-1}+\theta_\ell+\d\}_{i=1}^{ \ell-1}$}& $\supp(\theta_1)\sub\{\d_p\mid 1\leq p\leq n\}$ \\
&&\\
&&$\supp(\theta_\ell)\sub\{\d_p\mid 1\leq p\leq n\}$\\
 &&\\\cdashline{2-3}
 &\multirow{4}{*}{$\{-2\theta_1+\d,\theta_i-\theta_{i+1},\theta_{\ell-1}+\theta_\ell\}_{i=1}^{\ell-1}$}& $\supp(\theta_1)\sub\{\ep_i\mid 1\leq i\leq m\}$ \\
&&\\
&&$\supp(\theta_\ell)\sub\{\ep_i\mid 1\leq i\leq m\}$\\
 &&\\\cdashline{2-3}
 &\multirow{4}{*}{$\{-2\theta_1,\theta_i-\theta_{i+1},2\theta_\ell+\d\}_{i=1}^{\ell-1}$}& $\supp(\theta_1)\sub\{\d_p\mid 1\leq p\leq n\}$ \\
&&\\
&&$\supp(\theta_\ell)\sub\{\ep_i\mid 1\leq i\leq m\}$\\
 &&\\
\hline
\end{tabular}
\end{table}}

\bigskip

 Then $\Pi$ is a base of $R$ and each base of $R$ is of the form $\pm\Pi.$ In particular, each base of $R$ is an S-base. Moreover, for each type, bases in different rows are not conjugate under quasi-Weyl group $W$ of $R$ and  bases of the same form of each row are conjugate under $W.$
\end{Thm}

\section{Bases of twisted affine root supersystems}
Suppose $\gG$ is a twisted affine Lie superalgebra with Cartan subalgebra $\hH$ and root system $R;$ keep the same notation as in Introduction, recall  the form $\fm,$ $\ep_i$'s, $\d_p$'s,  Table~\ref{table1} and that $V=\sspan_\bbbr R$.
\begin{deft}
{\rm
 \begin{itemize}
 \item[(i)]
 A linearly independent subset of $R$ is called a {\it base} of $R$ if each element of $R$  can be written as a linear combination of elements of $\Pi$ with integeral coefficients such that all coefficients are nonnegative or all coefficients are nonpositive.
    The set of nonzero roots of $R\cap \sspan_{\bbbz^{\geq 0}}\Pi$   (resp.~$R\cap~\sspan_{\bbbz^{\leq 0}}\Pi$) is denoted by $R^+(\Pi)$ (resp. $R^-(\Pi)$) and  their elements are   called {\it positive} (resp. {\it negative}) roots.
 \item[(ii)]  A linearly independent subset $\Pi$ of $R$ is called an S-base if for each $\a\in \Pi,$ there are nonzero root vectors $x_\a$ and $y_\a$ corresponding to $\a$ and $-\a$ respectively such that $\{x_\a,y_\a\mid \a\in \Pi\}\cup \hH$ generates the affine Lie superalgebra $\gG$ and $[x_\a,y_\b]=0$ for $\a\neq \b.$
 \end{itemize}
    }
\end{deft}

\begin{lem}\label{reduce}
Suppose that    $\Pi$ is a base of $R.$ Then $\Pi\cap 2R=\Pi\cap \bbbz\d=\emptyset.$
\end{lem}
\pf If $\a,2\a\in R$ and $2\a\in \Pi,$ then $\a\in \frac{1}{2}\Pi$ which is  a contradiction. Also if $k\d\in \Pi$ for some $k\in\bbbz,$ then  for $\a\in \Pi\setminus \bbbz\d,$ $\a-4k\d\in R$ has been written as a $\bbbz$-linear combination of the elements of $\Pi$ with the opposite signs which is  a contradiction.\qed

\begin{lem}\label{new}
\begin{itemize}
\item[(i)] Each S-base of $R$ is a base of $R.$
\item[(ii)] Suppose that   $\sg$ is an even reflection. If  $\Pi$ is an S-base of $R$, then $\sg(\Pi)$ is also an S-base of $R$.
\item[(iii)] Suppose that $\Pi$ is an S-base of $R$ and $\sg$ is an odd reflection based on a nonsingular root of $\Pi.$ Then $\sg(\Pi)$ is also an S-base of $R$.
\item[(iv)]  Suppose that   $\Pi$  is an S-base  of $R.$ If $B$ is a base of $R$ such that either $R^+(\Pi)\cap R^-(B)$ or $R^+(\Pi)\cap R^+(B)$ is a finite set, then $B$ is an S-base.
\end{itemize}
\end{lem}
\pf (i) follows from the definition.

(ii) Suppose that $\a$ is a nonzero real root. Replacing $\a$ with $2\a$ if it is necessary, we assume $\a$ is an even real root. Fix $e_\a\in \gG^\a$ and $f_\a\in\gG^{-\a}$ such that $(e_\a,[e_\a,f_\a],f_\a)$ is an $\frak{sl}_2$-triple. Define
\[\theta_\a:=\hbox{exp}(\ad_{e_\a}) \hbox{exp}(-\ad_{f_\a}) \hbox{exp}(\ad_{e_\a}).\] Then $\theta_\a$ ia an automorphism of $\gG$ and
we have
\begin{equation}\label{theta}
\theta_\a(\hh)=\hh\andd \theta_\a(\gG^\b)=\gG^{r_\a(\b)}\quad\quad (\b\in R^\times_{re}\cup R_{ns}^\times ).
\end{equation}
Assume $\Pi$ is an S-base. Using part (i) and Lemma \ref{reduce}, we have $\Pi\sub R_{re}^\times\cup R_{ns}^\times.$ Fix weight vectors $x_\b\in\gG^\b$ and $y_\b\in\gG^{-\b}$ ($\b\in \Pi$) such that
$[x_\b,y_\gamma]=0$ for $\b\neq \gamma$ and $\{x_\b,y_\b\mid \b\in \Pi\}\cup \hh$ generates $\gG.$ Then we have $[\theta_\a(x_\b),\theta_\a(y_\gamma)]=0$ for $\b\neq \gamma$ and (\ref{theta}) implies that $\{\theta_\a(x_\b),\theta_\a(y_\b)\mid \b\in \Pi\}\cup \hh$ is a generating set for $\gG$. Since $\theta_\a(x_\b)\in \gG^{r_\a(\b)}$ and $\theta_\a(y_\b)\in \gG^{-r_\a(\b)},$ we are done.

(iii) See  Lemma  3.1 of \cite{serg2}.

(iv) Using Lemma~4.4  of \cite{serg2} together with parts (ii) and (iii) and the  same argument as in \cite[Lem.~4.5]{serg2}, we get that  $ B$ or $-B$ is obtained from $\Pi$ by even and odd reflections. Therefore parts (ii) and (iii) imply that $B$ is an S-base.

\medskip

We mention that the twisted affine root supersystem  $R$ is not necessarily  preserved by the quasi-Weyl group   of $R;$ e.g.,
in each of the types listed in Table \ref{table1}, $r_{\ep_i-\d_p}(2\d_p)=2\ep_i\not\in R,$ for $1\leq i\leq m$ and $1\leq p\leq n.$
But one can easily verify the following lemma.

\begin{lem}
\label{base}
Suppose that $w$ is an element of the quasi-Weyl group of $R$ with $w^{-1}(R)\sub R.$ If $\Pi$ is a base of $R,$ then $w(\Pi)$ is also a  base of $R.$
\end{lem}
\begin{cor}\label{belong-to}
 Assume $R=A(2k-1,2\ell-1)^{(2)},\;(k,\ell)\neq (1,1)$ and suppose that $p,q\in\bbbz$ and  $\zeta,\eta\in \{\pm\ep_i,\pm\d_j\mid 1\leq i\leq m,1\leq j\leq n\}$ with $\zeta\neq \pm\eta.$ Define
\begin{align*}
\mathscr{I}_{\zeta,p}:=&\left\{
\begin{array}{ll}
r_{2\zeta+p\d}& \hbox{if $\zeta\in\{\pm\d_j\mid 1\leq j\leq n\}$ and $p\in2\bbbz\d,$ }\\
r_{2\zeta+p\d} &  \hbox{if $\zeta\in\{\pm\ep_i\mid 1\leq i\leq  m\}$ and $p\in2\bbbz\d+\d,$}\\
%r_{2\zeta+(p-1)\d}& \hbox{if $p\in2\bbbz\d+\d$ and $\zeta\in\{\pm\d_p\mid 1\leq j\leq m\}$}\\
%r_{2\zeta+(p-1)\d} &  \hbox{if $p\in2\bbbz\d$ and $\zeta\in\{\pm\ep_i\mid 1\leq i\leq  n\}$}
\end{array}
\right.\\
\mathscr{J}_{\zeta,\eta,p}:=&r_{\zeta-\eta+p\d}\;\;\;\hbox{if $\zeta,\eta\in \{\pm\ep_i\mid 1\leq i\leq m\}$ or $\zeta,\eta\in \{\pm\d_j\mid 1\leq j\leq n\}$},\\
 \mathscr{S}_{\zeta,\eta,p,q}:=&r_{\zeta-\eta+p\d}r_{\zeta-\eta+q\d},\quad\andd
\mathscr{T}_{\zeta,\eta,p,q}:=r_{\zeta+\eta+p\d}r_{\zeta-\eta+q\d}.
\end{align*}
Then for $w\in \{\mathscr{I}_{2\zeta,p},\mathscr{S}_{\zeta,\eta,p,q},\mathscr{T}_{\zeta,\eta,p,q},\mathscr{J}_{\zeta,\eta,p}\},$ $w(R)\cup w^{-1}(R)\sub R;$ in particular,  if $\Pi$ is a base of $R,$ then $w(\Pi)$ is also a base of $R.$
\end{cor}
\pf %
%Then for  $r:=\mathscr{S}_{\zeta,\eta,k}(R), \mathscr{T}_{\zeta,\eta,k},$ $r(R)\sub R;$ in particular if $\Pi$  is  a base of $R,$ then $r(\Pi)$ is a base of $R$ as well.
%
%\[\mathscr{S}_{\zeta,\eta,k}(2\zeta+p\d)=r_{\zeta-\eta}(2\eta-2k\d+p\d)=2\zeta-2k\d+p\d\]
%
%\[\mathscr{T}_{\zeta,\eta,k}(2\zeta+p\d)=r_{\zeta+\eta}(2\eta-2k\d+p\d)=-2\zeta-2k\d+p\d\]
%
%we do not need to check short roots.
It easily follows using a direct verification together with Lemma \ref{base}.\qed

\bigskip

\subsection{Type $\boldsymbol{A(2m,2n)^{(4)}}$}
Set $$\dot R:=
  \pm\{\ep_i,\d_j,2\ep_i,2\d_j,\ep_i\pm\ep_r,\d_j\pm\d_s,\d_j\pm\ep_i\mid 1\leq i\neq r\leq m,1\leq j\neq s\leq n\},
$$
$$\begin{array}{ll}
\dot R_{sh}:= \{\pm\ep_i,\pm\d_p\mid 1\leq i\leq m,\;1\leq p\leq n\}&(\hbox{short roots}),\\
\dot R_{lg}:=\{\pm\ep_i\pm\ep_j,\pm\d_p\pm\d_q\mid 1\leq i\neq j\leq m,\;1\leq p\neq q\leq n\}&(\hbox{long roots}),\\
\dot R_{ex}:=\{\pm2\ep_i,\pm2\d_p\mid 1\leq i\leq m,\;1\leq p\leq n\}&(\hbox{extra long roots}),\\
\dot R_{ns}:=\{\pm\ep_i\pm\d_p\mid1\leq i\leq m,\;1\leq p\leq n\}
\end{array}
$$ and   $$R:=\bbbz\d\cup  (\dot R_{sh}+\bbbz\d)\cup (\dot R_{lg}\cup \dot R_{ns})+2\bbbz\d)\cup (\{\pm2\ep_i\}_{i=1}^m+4\bbbz\d+2\d)\cup (\{\pm2\d_p\}_{p=1}^n+4\bbbz\d).$$
For $$S:=\bbbz\d\cup  (\dot R_{sh}+\bbbz\d)\cup ((\dot R_{lg}\cup \dot R_{ns})+2\bbbz\d),$$ $(\sspan_\bbbr S,\fm_*,S)$ is the root system of $D^{(2)}_{m+n+1}.$  Using Lemma \ref{reduce}, each base of $R$ is a base of $S$ and conversely each base of $S$ is a base of $R.$
So for $\theta_i$ ($1\leq i\leq \ell$) as in Theorem \ref{main} (see Table \ref{table2} for further assumptions), we have
\[S=\bbbz\d\cup (\{\pm\theta_i\mid 1\leq i\leq \ell\}+\bbbz\d)\cup (\{\pm\theta_i\pm\theta_j\mid 1\leq i\neq j\leq \ell\}+2\bbbz\d)\] and
\[\Pi=\{\theta_i-\theta_{i+1},\theta_\ell,\d-\theta_1,\mid 1\leq i\leq \ell-1\}\] is a base of $S$ (and so a base of $R$). Moreover, the quasi-Weyl group of $R$ is the Weyl group of $S$ and so from affine Lie theory, up to $\pm1$-multiple, all bases of $R$  are conjugate under quasi-Weyl group; in particular, all bases of $R$ are of the form $\Pi$ or $-\Pi.$  Next suppose $\D$ is the standard base of $R$ (which is also a  base of $S$), then  again using affine Lie theory, for each base $\D'$ of $S$ (or equivalently of $R$), we have $S^+(\D)\cap S^-(\D')$ or  $S^+(\D)\cap S^+(\D')$ is finite and consequently either   $R^+(\D)\cap R^-(\D')$  or $R^+(\D)\cap R^+(\D')$  is finite. Therefore,  we get using Lemma \ref{new} that  $\D$ is an S-base. This in particular completes the proof of  Theorem \ref{main} for $R=A(2m,2n)^{(4)}$.
\subsection{Type  $\boldsymbol{D(m+1,n)^{(2)}}$} One can prove Theorem \ref{main} using the same argument as in the previous case.
\begin{comment}
As before assume $m=0$ if  $R$ is of type $C(n+1)^{(2)}$ and set $$\hbox{\small $R:= \bbbz\d\cup  \bbbz\d\pm\{\ep_i,\d_j\mid 1\leq i\leq m,1\leq j\leq n\}\cup 2\bbbz\d\pm\{2\d_j,\ep_i\pm\ep_r,\d_j\pm\d_s,\d_j\pm\ep_i\mid i\neq r,j\neq s\}$}$$ and $$S:=\bbbz\d\cup  \bbbz\d\pm\{\ep_i,\d_j\mid 1\leq i\leq m,1\leq j\leq n\}\cup 2\bbbz\d\pm\{\ep_i\pm\ep_r,\d_j\pm\d_s,\d_j\pm\ep_i\mid i\neq r,j\neq s\}.$$ As in the previous case,  $(\sspan_\bbbr S,\fm_*,S)$ is an affine root system of type $B_{m+n}$ and using Lemma \ref{reduce}, each base of $R$ is a base of $S$ and conversely each base of $S$ is a base of $R;$ in particular
\[\Pi:=\{\ep_i-\ep_{i+1},\ep_{m}-\d_1,\d_p-\d_{p+1},-2\ep_1+\d\mid 1\leq i\leq m,1\leq p\leq n\}\] is a base of $R$ and each base of $R$ is conjugate with $\pm\Pi$ under quasi-Weyl group of $R.$

So if
$k_1,\ldots,k_\ell\in\bbbz,$ $\zeta_1,\ldots\zeta_\ell\in\{\pm \ep_i,\pm\d_p\mid 1\leq i\leq m,1\leq p\leq n\}$ with $\zeta_i\neq \zeta_j$ for $i\neq j$ such that for $\theta_i:=\zeta_i+k_i\d$ ($1\leq i\leq \ell$),
\[B=\{\theta_i-\theta_{i+1},\theta_\ell,\d-2\theta_1,\mid 1\leq i\leq \ell-1\}\] is a base of $R$ and each base of $R$ is of the form $B$ or $-B$ and we have   \[R^+(B)=S^+(B)\cup (R\cap 2S^+(B)).\] If $\D$ and $\D'$ are two bases of the form $B,$ then $R^+(\D)\cap R^-(\D)$ is a finite set so if $\D$ is the standard base, using \cite[Cor. and Lem.~3.1]{serg2}, we get $\D'$  is an S-base.
\end{comment}
\subsection{Type $\boldsymbol{A(2m,2n-1)^{(2)}}$}
Set $$\hbox{\small$\dot R:=\{0,\pm\ep_i,\pm\ep_i\pm\ep_j,\pm2\ep_i,\pm\d_p,\pm\d_p\pm\d_q,\pm2\d_p,\pm\ep_i\pm\d_p\mid 1\leq i\neq j\leq m,\;1\leq p\neq q\leq n\}$},$$
$$\begin{array}{ll}
\dot R_{sh}:= \{\pm\ep_i,\pm\d_p\mid 1\leq i\leq m,\;1\leq p\leq n\}&(\hbox{short roots}),\\
\dot R_{lg}:=\{\pm\ep_i\pm\ep_j,\pm\d_p\pm\d_q\mid 1\leq i\neq j\leq m,\;1\leq p\neq q\leq n\}&(\hbox{long roots}),\\
\dot R_{ex}:=\{\pm2\ep_i,\pm2\d_p\mid 1\leq i\leq m,\;1\leq p\leq n\}&(\hbox{extra-long roots}),\\
\dot R_{ns}:=\{\pm\ep_i\pm\d_p,1\leq i\leq m,\;1\leq p\leq n\}.
\end{array}
$$ and   $$R:=\bbbz\d\cup  ((\dot R\setminus\dot R_{ex})+\bbbz\d)\cup (\{\pm2\ep_i\}_{i=1}^m+2\bbbz\d+\d)\cup (\{\pm2\d_p\}_{p=1}^n+2\bbbz\d).$$ We recall that
\[\ell=m+n\]  and set
\begin{align} \label{t}
& \dot S:=\dot R\setminus\dot R_{ex}=\{0,\pm\zeta,\pm\zeta\pm\xi\mid \zeta\neq \xi\in\{\ep_i,\d_p\mid 1\leq i\leq m,\;1\leq p\leq n\}\},\nonumber \\
 & S:=\dot S+\bbbz\d\andd\\
 & T:=S\cup (\{\pm2\zeta\mid \zeta\in\{\ep_i,\d_p\mid 1\leq i\leq m,\;1\leq p\leq n\}\}+2\bbbz\d+\d).\nonumber
\end{align}
Then  $(V,\fm_{*},S)$ is the root system of $B_{m+n}^{(1)}$ and
$(V,\fm_{*},T)$ is the  root system of $A_{2(m+n)}^{(2)}.$   We have
 \begin{align}
\hbox{\small $T\setminus R=\{\pm2\d_p\mid 1\leq p\leq n\}+2\bbbz\d+\d\andd R\setminus T=\{\pm2\d_p\mid 1\leq p\leq n\}+2\bbbz\d.$}\label{t-r}
\end{align}

In what follows, we call elements of $(\{\pm2\ep_i\}_{i=1}^n+2\bbbz\d)\cup (\{\pm2\d_p\}_{p=1}^m+2\bbbz\d+\d),$ extra-long roots of $R$ and denote by   $W$  the quasi-Weyl  group of $R.$
\begin{Pro}\label{with-long}
If $\Pi\sub R$ is a base of $T,$ then $\Pi$ is a base of $R$ and conversely, each base of  $R$ containing an extra long root is a base of $T.$
\end{Pro}
\pf  Suppose that  $\Pi\sub R$ is a base of $T.$   Since for each $1\leq p\leq n,$ $2\d_p+2\bbbz\d=2(\d_p+\bbbz\d)\sub  2T,$ $\Pi$ is  a base of $R.$

Conversely, suppose that  $\Pi$ is a base of $R$ containing an extra long root $\a,$ then by Lemma \ref{reduce},  $\a=2r\ep_i+2k\d+\d$ for some  $r\in\{\pm1\},$ $1\leq i\leq m$ and $k\in\bbbz.$ Set $\dot\ep_i:=-r\ep_i-k\d\in R,$ so $$\a=-2\dot\ep_i+\d\in \Pi.$$ To the contrary assume  $\Pi$ is not a base of $T.$
Using (\ref{t-r}), there is an index $p_0,$ $s\in\{\pm1\}$ and $k'\in\bbbz$  such that $2s\d_{p_0}+(2k'\d+\d)\not\in \hbox{span}_{\bbbz^{\geq0}}\Pi\cup \hbox{span}_{\bbbz^{\leq0}}\Pi.$ Set $\dot\d_{p_0}:=s\d_{p_0}+k'\d\in R.$ So
\begin{equation}
\label{i0}
2\dot\d_{p_0}+\d\not\in \hbox{span}_{\bbbz^{\geq0}}\Pi\cup \hbox{span}_{\bbbz^{\leq0}}\Pi.
\end{equation}
 Since $\Pi$ is a basis of $V,$ by Lemma \ref{reduce}, there is an element of $\Pi$ of the form $t\dot\d_{p_0}+\rho$ for some $t\in\{\pm1\}$ and  $\rho\in\bbbz\d\cup(\{\pm\d_q\mid 1\leq q\neq p_0\leq n\}+\bbbz\d)\cup(\{\pm\ep_j\mid 1\leq j\leq n\}+\bbbz\d).$
 Assume $\Pi=\{\a_1,\ldots,\a_{\ell+1}\}$ with $$\a_1=-2\dot\ep_i+\d\andd \a_2=t\dot\d_{p_0}+\rho.$$
We know that $\dot\d_{p_0}+\dot\ep_i,t\dot\d_{p_0}-\rho+t\d \in R.$ Therefore, there are integers $r_i,s_i$ ($1\leq i\leq \ell+1$) such that
$$\dot\d_{p_0}+\dot \ep_i=\sum_{i=1}^{\ell+1}r_i\a_i\andd t\dot\d_{p_0}-\rho+t\d=\sum_{i=1}^{\ell+1}s_i\a_i$$
and $r_i$'s (resp. $s_i$'s) are all non-negative or all non-positive.
We have
\begin{align*}
2\dot\d_{p_0}+\d=(-2\dot\ep_i+\d)+2(\dot\d_{p_0}+\dot\ep_i)
=\a_1+\sum_{i=1}^{\ell+1}2r_i\a_i= (1+2r_1)\a_1+\sum_{i=2}^{\ell+1}2r_i\a_i
\end{align*}
and
 \begin{align*}
2t\dot\d_{p_0}+t\d=(t\dot\d_{p_0}+\rho)+(t\dot\d_{p_0}-\rho+t\d)=\a_2+\sum_{i=1}^{\ell+1}s_i\a_i
=(1+s_2)\a_2+\sum_{2\neq i=1}^{\ell+1}s_i\a_i.
\end{align*}
If $t\dot\d_{p_0}-\rho+t\d$ is a positive root that is  for all $i,$ $s_i\geq 0,$ then $2t\dot\d_{p_0}+t\d\in \hbox{span}_{\bbbz^{\geq0}}\Pi$ contradicting (\ref{i0}). So $t\dot\d_{p_0}-\rho+t\d$  is a negative root; in other words $s_i\leq 0$ for all $i.$ Using (\ref{i0}), we get  that $0< 1+s_2\leq 1.$ This implies that $s_2=0.$ Therefore, $$(1+2r_1)\a_1+\sum_{i=2}^{\ell+1}2r_i\a_i=2\dot\d_{p_0}+\d=t\a_2+\sum_{2\neq i=1}^{\ell+1}ts_i\a_i$$ which implies that $2r_2=t\in\{\pm1\},$ a contradiction.\qed

\begin{lem}\label{pos-neg}
Assume  $\theta_i$'s and $\Pi$ are as in  Theorem \ref{main}. Then $\Pi$ is a base of $R$ and  \begin{align*}
R^+(\Pi)=&\{\pm2\theta_i\mid \supp(\theta_i)\sub\{\ep_j\mid 1\leq j\leq m\}\}+2\bbbz^{\geq0}\d+\d\\
\cup & \{\pm2\theta_i\mid \supp(\theta_i)\sub\{\d_p\mid 1\leq p\leq n\}+2\bbbz^{>0}\d\\
\cup &\{0,\pm\theta_i,\pm\theta_i\pm\theta_j\mid 1\leq i<j\leq \ell\}+\bbbz^{>0}\d\\
\cup &\{\theta_i,\theta_i\pm \theta_j\mid 1\leq i<j\leq \ell\}\cup \{2\theta_i\mid \supp(\theta_i)\sub\{\d_p\mid 1\leq p\leq n\}\}.
\end{align*}
\end{lem}
\pf
We have two cases:

\noindent{Case 1.} $\Pi=\{ \d-2\theta_1,\theta_i-\theta_{i+1},\theta_\ell\mid 1\leq i\leq \ell-1\}:$ In this case, the result follows  easily from Proposition \ref{with-long} and affine Lie theory as $\Pi$ is a  base of
\[
T=\{0,\pm\theta_i,\pm\theta_i\pm\theta_j\mid 1\leq i\neq j\leq \ell\}+\bbbz\d\cup\{\pm2\theta_i\mid 1\leq i\leq \ell\}+2\bbbz\d+\d,\] which is the root system of $A_{2\ell}^{(2)}.$

\noindent{Case 2.} $\Pi=\{\d-(\theta_1+\theta_2),\theta_i-\theta_{i+1},\theta_\ell\mid 1\leq i\leq \ell-1\}:$ We have
\begin{align*}
S=& \{0,\pm\ep_i,\pm\d_p,\pm\ep_i\pm\ep_j,\pm\d_p\pm\d_q,\pm\ep_i\pm\d_p\mid 1\leq i\neq j\leq m,\;1\leq p\neq q\leq n\}+\bbbz\d\\
=&\{0,\pm\theta_i,\pm\theta_i\pm\theta_j\mid 1\leq i\neq j\leq \ell\}+\bbbz\d
\end{align*}
 which is the root system of $B_\ell^{(1)}$ and that  $\Pi$ is a base of $S$ with
\[S^+(\Pi)=\{0,\pm\theta_i,\pm\theta_i\pm\theta_j\mid 1\leq i<j\leq \ell\}+\bbbz^{>0}\d\cup \{\theta_i,\theta_i\pm\theta_j\mid 1\leq i<j\leq \ell\}.\]
 So we have
\begin{equation}\label{com}-2\theta_2+\d=\theta_1-\theta_2+(-(\theta_1+\theta_2)+\d)\in \Pi+\Pi\sub \hbox{span}_{\bbbz^{\geq 0}}\Pi.
\end{equation}
Therefore,  for $t\in\{\pm 1\},$ $k\in\bbbz^{\geq0}$ and $2\leq i\leq m+n,$ we have
\begin{align*}
2t\theta_i+2k\d+\d=2(t\theta_i+\theta_2)+(-2\theta_2+\d)+2k\d\in \hbox{span}_{\bbbz^{\geq 0}}\Pi.
\end{align*}
We are done using this together with the fact that $\supp(\theta_1)\sub\{\d_p\mid 1\leq p\leq n\}.$
\qed

\medskip

The following theorem  completes the proof of Theorem \ref{main} for the case under consideration:
\begin{Thm}
\begin{itemize}
\item[(i)]
Each base of $R$ is of the form $\pm\Pi$ where $\Pi$ is as in Table \ref{table2}. In particular, if
\begin{equation}\label{ep-delta}\xi_i:=\left\{\begin{array}{ll}
\d_i& 1\leq i\leq n,\\
\ep_{i-n}& n+1\leq i\leq \ell,
\end{array}\right.
\;\; \eta_i:=
\left\{\begin{array}{ll}
\ep_i& 1\leq i\leq m,\\
\d_{i-m}& m+1\leq i\leq \ell,
\end{array}\right.
\end{equation}then
\[
\Pi_1:=\{\d-2\eta_1,\eta_i-\eta_{i+1},\eta_\ell\}\andd
\Pi_2:=\{\d-(\xi_1+\xi_2),\xi_i-\xi_{i+1},\xi_{\ell}\}
\]
are bases of $R$ and up to conjugacy under quasi-Weyl group, $\pm\Pi_1$ and $\pm\Pi_2$ are the only bases of $R.$
\item[(ii)] Each base of $R$ is an S-base.
\end{itemize}
\end{Thm}
\pf (i) Suppose that $B$ is a base of $R.$ Two cases can happen: Either $B$ contains an extra long root  or it  is  contained in $S$.

\noindent{ Case 1.} \underline{ $B$ contains an extra long root}: We first mention that  $\Pi_1$ is a base of $T$ and that under the Weyl group action of $T,$ each base of $T$ is conjugate with $\Pi_1$ or $-\Pi_1.$ If $W_1$ is the Weyl group of $T,$ then  $W_1$ equals with the group generated by the reflections  based on the elements of $\Pi_1.$ So $W_1\sub W.$

By Proposition \ref{with-long}, B is a base of $T;$ in particular,
there is an element $w\in W_1\sub W$ such that either $w(\Pi_1)=B$ or $w(\Pi_1)=-B.$ Set
\[\theta_i:=w(\eta_i)\in\{\pm \ep_j,\pm\d_p\mid 1\leq j\leq m,1\leq p\leq n\}+\bbbz\d\quad\quad(1\leq i\leq \ell).\]
Since $-2\theta_1+\d=w(-2\eta_1+\d)\in \pm B\sub R,$ we get that $\theta_1\in \{\pm \ep_i\mid 1\leq i\leq m\}+\bbbz\d.$ Also we have
$$\pm B=w(\Pi_1)=\{\theta_1-\theta_2,\ldots,\theta_{\ell-1}-\theta_{\ell},\theta_{\ell},\d-2\theta_1\}.$$

\noindent{Case 2.} \underline{$B$ is a base of $R$ contained in $S$}: In this case,  $B$ is a base of $S;$ in particular,  $B$ is either conjugate with $\Pi_2$ or $-\Pi_2$ under the  Weyl group action  of $S.$ But Weyl group of $S$ is a subgroup of $W,$ so there is $w\in W$ with $w(\Pi_2)=\pm B.$
Set $$\theta_i:=w(\eta_i)\in \{\pm\ep_t,\pm\d_p\mid 1\leq t\leq m,1\leq p\leq n\}+\bbbz\d\;\;\;\;(1\leq i\leq \ell).$$
 So
$$\pm B=w(\Pi_2)=\{\theta_1-\theta_2,\ldots,\theta_{\ell-1}-\theta_{\ell},\theta_{\ell},\d-(\theta_{1}+\theta_2)\}.$$
We have
$$-2\theta_1+\d=-(\theta_1-\theta_2)+(\d-(\theta_1+\theta_2))\in-\Pi_2+\Pi_2.$$ This means that
  $-2\theta_1+\d \not \in R,$ i.e., $\theta_1\in \{\pm\d_p\mid 1\leq p\leq n\}+\bbbz\d.$ This completes the proof.

(ii) Suppose that $\D$ is the standard base of $R$, then by part (i) and  Lemma \ref{pos-neg}, for each base $\D'$, either $R^+(\D)\cap R^-(\D')$ or $R^+(\D)\cap R^+(\D')$  is finite. Therefore, we get the result using Lemma \ref{new}.\qed

\subsection{\bf Type $\boldsymbol{A(2m-1,2n-1)^{(2)}}$}
Suppose that $m,n$ are two positive integers with $(m,n)\neq (1,1)$
and  assume $R$ is the  root system of  $A(2m-1,2n-1)^{(2)};$ in other words, \[\hbox{\small$R=(\{0,\pm\ep_i\pm\ep_j,\pm\d_p\pm\d_q,\pm\ep_i\pm\d_p\}_{i,p}+\bbbz\d)\cup (\{\pm2\ep_i\}_i+2\bbbz\d+\d)\cup (\{\pm2\d_p\}_p+2\bbbz\d)$}\]
where $i\neq j$ run over $\{1,\ldots,m\}$ and $p\neq q$ run over $\{1,\ldots,n\}.$ Set
\begin{equation}
\label{long}R_{lg}:=(\{\pm2\ep_i\}_{i=1}^m+2\bbbz\d+\d)\cup (\{\pm2\d_p\}_{p=1}^n+2\bbbz\d)\quad\quad(\hbox{long roots})
\end{equation}
 and recall   $\fm_*$ from (\ref{form*}) in  Introduction  and that
\[\ell=m+n.\]
\begin{Pro}\label{bases-last}
Recall  $\theta_i$'s and $\Pi$  from Theorem \ref{main}.
\begin{itemize}
\item[(i)] If $\Pi=\{\d-(\theta_1+\theta_2),\theta_i-\theta_{i+1},\theta_{\ell-1}+\theta_\ell\mid 1\leq i\leq \ell-1\},$ then $\Pi$ is a base of $R$ and
    \begin{align*}
    R^+(\Pi)=&\{\theta_i\pm\theta_j\mid i<j\}\cup \{2\theta_i\mid \supp(\theta_i)\sub\{\d_p\mid 1\leq j\leq n\}\}\\
    \cup& \{0,\pm\theta_i\pm\theta_j\mid i\neq j\}+\bbbz^{>0}\d \\
    \cup& \{\pm 2\theta_i\mid \supp(\theta_i)\sub\{\ep_j\mid 1\leq j\leq m\}\}+2\bbbz^{\geq 0}\d+\d\\
    \cup &\{\pm2\theta_i\mid \supp(\theta_i)\sub\{\d_p\mid 1\leq j\leq n\}\}+2\bbbz^{>0}\d.
    \end{align*}
\item[(ii)] If    $\Pi=\{-2\theta_1,\theta_i-\theta_{i+1},\theta_{\ell-1}+\theta_\ell+\d\mid 1\leq i\leq \ell-1\},$ then $\Pi$ is  a base of $R$ and
\begin{align*}
    R^+(\Pi)=&\{\theta_i-\theta_j,-\theta_i-\theta_j\mid i<j\}\cup \{-2\theta_i\mid \supp(\theta_i)\sub\{\d_p\mid 1\leq j\leq n\}\}\\
    \cup& \{0,\pm\theta_i\pm\theta_j\mid i\neq j\}+\bbbz^{>0}\d \\
    \cup& \{\pm 2\theta_i\mid \supp(\theta_i)\sub\{\ep_j\mid 1\leq j\leq m\}\}+2\bbbz^{\geq 0}\d+\d\\
    \cup &\{\pm2\theta_i\mid \supp(\theta_i)\sub\{\d_p\mid 1\leq j\leq n\}\}+2\bbbz^{>0}\d.
    \end{align*}
    \item[(iii)] If $\Pi=\{-2\theta_1+\d,\theta_i-\theta_{i+1},\theta_{\ell-1}+\theta_\ell\mid 1\leq i\leq \ell-1\},$ then $\Pi$ is  a base of $R$ and
\begin{align*}
    R^+(\Pi)=&\{\theta_i\pm\theta_j\mid i<j\}\cup \{2\theta_i\mid \supp(\theta_i)\sub\{\d_p\mid 1\leq j\leq n\}\}\\
    \cup& \{0,\pm\theta_i\pm\theta_j\mid i\neq j\}+\bbbz^{>0}\d \\
    \cup& \{\pm 2\theta_i\mid \supp(\theta_i)\sub\{\ep_j\mid 1\leq j\leq m\}\}+2\bbbz^{\geq 0}\d+\d\\
    \cup &\{\pm2\theta_i\mid \supp(\theta_i)\sub\{\d_p\mid 1\leq j\leq n\}\}+2\bbbz^{>0}\d.
    \end{align*}
\item[(iv)] If  $\Pi=\{-2\theta_1,\theta_i-\theta_{i+1},2\theta_\ell+\d\mid 1\leq i\leq \ell-1\},$ then $\Pi$ is  a base of $R$ and
\begin{align*}
    R^+(\Pi)=&\{\theta_i-\theta_j,-\theta_i-\theta_j\mid i<j\}\cup \{-2\theta_i\mid \supp(\theta_i)\sub\{\d_p\mid 1\leq j\leq n\}\}\\
    \cup& \{0,\pm\theta_i\pm\theta_j\mid i\neq j\}+\bbbz^{>0}\d \\
    \cup& \{\pm 2\theta_i\mid \supp(\theta_i)\sub\{\ep_j\mid 1\leq j\leq m\}\}+2\bbbz^{\geq 0}\d+\d\\
    \cup &\{\pm2\theta_i\mid \supp(\theta_i)\sub\{\d_p\mid 1\leq j\leq n\}\}+2\bbbz^{>0}\d.
    \end{align*}
\end{itemize}
\end{Pro}
\pf (i) Setting
\begin{align}
S:=&\pm\{0,\theta_i\pm\theta_j\mid 1\leq i\neq j\leq \ell\}+\bbbz\d\label{S}\\
=&\pm\{0,\ep_i\pm\ep_j,\d_p\pm\d_q,\ep_i\pm\d_p\mid 1\leq i\neq j\leq m,1\leq p\neq q\leq n\}+\bbbz\d,\nonumber
\end{align}
we get that $(\sspan_\bbbr S,\fm_*,S)$ is the root system of $D^{(1)}_\ell$ and $\Pi$ is a base of $S.$ One knows that the set  $S^+$ of positive roots of $S$ with respect to $\Pi$ is
\begin{equation}\label{plus}
S^+=\{\theta_i\pm\theta_j\mid 1\leq i<j\leq \ell\}\cup (\pm\{0,\theta_i\pm\theta_j\mid1\leq  i\neq j\leq \ell\}+\bbbz^{>0}\d).
\end{equation} To complete the proof, it is enough  to show that  for all $1\leq i\leq \ell,$
\begin{itemize}
\item $\pm 2\theta_i+2\bbbz^{\geq 0}\d+\d\in \sspan_{\bbbz^{\geq0}}\Pi$ if $\supp(\theta_i)\sub\{\ep_j\mid 1\leq j\leq m\}$ and
\item $2\theta_i,\pm2\theta_i+2\bbbz^{>0}\d\in \sspan_{\bbbz^{\geq0}}\Pi$ if $\supp(\theta_i)\sub\{\d_p\mid 1\leq p\leq n\}.$
\end{itemize}
We first assume $\supp(\theta_i)\sub \{\ep_j\mid 1\leq j\leq m\}.$ Then $i\neq 1$ and so for $k\geq 0,$ using (\ref{plus}),  we have
\[\pm2\theta_i+(2k+1)\d=(-\theta_1\pm\theta_i+(2k+1)\d)+(\theta_1\pm\theta_i)\in \sspan_{\bbbz^{\geq 0}}\Pi.\]
Also if $\supp(\theta_i)\sub \{\d_p\mid 1\leq p\leq n\},$ we have $i\neq \ell$ and so   by (\ref{plus}), we have
\[\pm2\theta_i+2k\d=(\pm\theta_i+\theta_\ell+k\d)+(\pm\theta_i-\theta_\ell+k\d)\in \sspan_{\bbbz^{\geq 0}}\Pi\quad\quad (k> 0)\] and
\[2\theta_i=(\theta_i-\theta_\ell)+(\theta_i+\theta_\ell)\in \sspan_{\bbbz^{\geq 0}}\Pi\] and so we are done.

(ii)-(iv) Suppose $\Pi$ is as in the statements. To simplify our argument, we also denote the $\Pi$ introduced in (ii), (iii) and (iv) respectively by $B^2,$ $B^3$ and $B^4.$ Let $k\in\bbbz^{>0}$ and $1\leq i\neq j\leq \ell,$ then we have the following table:
{\footnotesize\[
\begin{tabular}{|c|l|c|}
\hline
$\theta_i-\theta_j$\;\;\footnotesize{$(i<j)$} &$(\theta_i-\theta_{i+1})+\cdots+(\theta_{j-1}-\theta_j)$ & $\in \sspan_{\bbbz^{\geq0}}\Pi$ \\
 \hline
\multirow{4}{*}{$\theta_i+\theta_j$\;\;\footnotesize{$(i<j)$}} &\multirow{2}{*}{ $(\theta_i-\theta_{\ell-1})+(\theta_{j}-\theta_{\ell})+(\theta_{\ell-1}+\theta_\ell)$} & \multirow{2}{*}{ $\in \sspan_{\bbbz^{\geq0}}B^3$}\\&&\\\cdashline{2-3}
&\multirow{2}{*}{ $-(-2\theta_1)-2(\theta_1-\theta_i)-(\theta_i-\theta_j)$} & \multirow{2}{*}{ $\in \sspan_{\bbbz^{\leq0}}B^2\cup \sspan_{\bbbz^{\leq0}}B^4$}\\&&\\
 \hline
\multirow{6}{*}{$\theta_i+\theta_j+\d$\;\;\footnotesize{$(i<j)$}}
&\multirow{2}{*}{ $(\theta_i-\theta_{\ell-1})+(\theta_j-\theta_\ell)+(\theta_{\ell-1}+\theta_\ell+\d)$} & \multirow{2}{*}{ $\in \sspan_{\bbbz^{\geq0}}B^2$}\\&&\\\cdashline{2-3}
&\multirow{2}{*}{ $(\theta_i+\theta_j)+(\theta_1+\theta_2)+(\theta_1-\theta_2)+(-2\theta_1+\d)$} & \multirow{2}{*}{ $\in \sspan_{\bbbz^{\geq0}}B^3$}\\&&\\\cdashline{2-3}
&\multirow{2}{*}{ $(\theta_i-\theta_{\ell})+(\theta_j-\theta_\ell)+(2\theta_\ell+\d)$} & \multirow{2}{*}{ $\in  \sspan_{\bbbz^{\geq0}}B^4$}\\
&&\\
 \hline
\multirow{8}{*}{$2\theta_i+\d$}
 &\multirow{2}{*}{$(\theta_i-\theta_{i+1})+(\theta_i+\theta_{i+1}+\d)$\;\;\footnotesize{(if $i\neq \ell)$}}
 & \multirow{2}{*}{$\in \sspan_{\bbbz^{\geq0}}B^2$} \\&& \\\cdashline{2-3}
 &\multirow{2}{*}{$2(\theta_i+\theta_{i+1})+2(\theta_1-\theta_{i+1})+(-2\theta_1+\d)$\;\;\footnotesize{(if $i\neq \ell)$}}&\multirow{4}{*}{ $\in \sspan_{\bbbz^{\geq0}}B^3$ }\\&&\\
&\multirow{2}{*}{$(-2\theta_1+\d)+2(\theta_1+\theta_\ell)$\;\;\footnotesize{(if $i= \ell)$}}
 &\\ &&\\\cdashline{2-3}
 &\multirow{2}{*}{$2(\theta_i-\theta_\ell)+(2\theta_\ell+\d)$}& \multirow{2}{*}{$\in \sspan_{\bbbz^{\geq0}}B^4$}\\&& \\\hline
\multirow{4}{*}{$-2\theta_i$} & \multirow{2}{*}{$-(\theta_i-\theta_{i+1})-(\theta_i+\theta_{i+1})$\;\;\footnotesize{(if $i\neq \ell)$}} & \multirow{2}{*}{$\in \sspan_{\bbbz^{\leq0}}B^3$}\\&&\\\cdashline{2-3}
& \multirow{2}{*}{$(-2\theta_1)+2(\theta_1-\theta_2)+\cdots+2(\theta_{i-1}-\theta_i)$} & \multirow{2}{*}{$\in \sspan_{\bbbz^{\geq0}}B^2\cup \sspan_{\bbbz^{\geq0}}B^4$}\\&&\\
\hline
\multirow{4}{*}{$\d$} & \multirow{2}{*}{$(\theta_1+\theta_2)+(\theta_1-\theta_2)+(-2\theta_1+\d)$} &\multirow{2}{*}{ $\in \sspan_{\bbbz^{\geq0}}B^3$}\\&&\\\cdashline{2-3}
& \multirow{2}{*}{$(2\theta_1+\d)+(-2\theta_1)$} &\multirow{2}{*}{ $\in \sspan_{\bbbz^{\geq0}}B^2\cup \sspan_{\bbbz^{\geq0}}B^4$}\\&&\\
\hline
\multirow{4}{*}{$\theta_i+\theta_j+k\d$} &\multirow{2}{*}{ $(\theta_i+\theta_j)+k\d$} & \multirow{2}{*}{ $\in \sspan_{\bbbz^{\geq0}}B^3$}\\&&\\\cdashline{2-3}
&\multirow{2}{*}{ $(\theta_i+\theta_j+\d)+(k-1)\d$} & \multirow{2}{*}{ $\in \sspan_{\bbbz^{\geq0}}B^2\cup \sspan_{\bbbz^{\geq0}}B^4$}\\&&\\
\hline
\multirow{4}{*}{ $-2\theta_i+k\d$} &\multirow{2}{*}{  $(-2\theta_1+\d)+2(\theta_1-\theta_i)+(k-1)\d$}&\multirow{2}{*}{  $\in \sspan_{\bbbz^{\geq0}}B^3$}\\&&\\\cdashline{2-3}
&\multirow{2}{*}{  $(-2\theta_i)+k\d$}&\multirow{2}{*}{  $\in \sspan_{\bbbz^{\geq0}}B^2\cup \sspan_{\bbbz^{\geq0}}B^4$}\\&&\\
 \hline
 \multirow{4}{*}{$\theta_i-\theta_j+k\d$} &\multirow{2}{*}{ $(\theta_i+\theta_j)+(-2\theta_j+\d)+(k-1)\d$} & \multirow{2}{*}{ $\in \sspan_{\bbbz^{\geq0}}B^3$}\\&&\\\cdashline{2-3}
&\multirow{2}{*}{ $(\theta_i+\theta_j+\d)+(-2\theta_j)+(k-1)\d$} & \multirow{2}{*}{ $\in \sspan_{\bbbz^{\geq0}}B^2\cup \sspan_{\bbbz^{\geq0}}B^4$}\\&&\\
 \hline
 \multirow{4}{*}{$2\theta_i+k\d$\;\;\footnotesize{($k\geq 2)$}}
  &\multirow{2}{*}{$(2\theta_i+\d)+(k-1)\d$}
 & \multirow{2}{*}{$\in \sspan_{\bbbz^{\geq0}}B^3\cup \sspan_{\bbbz^{\geq0}}B^4$} \\&& \\\cdashline{2-3}
 &\multirow{2}{*}{ $(\theta_i-\theta_1+(k-1)\d)+(\theta_i+\theta_1+\d)$}&\multirow{2}{*}{ $\in \sspan_{\bbbz^{\geq0}}B^2$ }\\&&\\\hline
 \multirow{2}{*}{$\theta_i+\theta_j-k\d$\;\;\footnotesize{($i<j)$}}
  &\multirow{2}{*}{$-(\theta_i-\theta_j)-(-2\theta_i+k\d)$}
 & \multirow{2}{*}{$\in \sspan_{\bbbz^{\leq0}}\Pi$} \\&&\\\hline
\end{tabular}
\]}

This completes the proof.\qed

We shall show that  the bases introduced in Proposition \ref{bases-last} are the only bases of $R.$
Let us start with stating two facts which  are easily verified and we use them   frequently in the sequel:
 \medskip

\noindent{\bf Fact 1.} If $\a,\b\in R$ with $(\a,\b)_*=0$ and $\supp(\a)\cap \supp(\b)\neq \emptyset,$ then $\supp(\a)= \supp(\b)$ and $|\supp(\a)|=|\supp(\b)|=2,$ in which $|\cdot|$ indicates the cardinal number; moreover, if $\supp(\a)=\{\zeta_1,\zeta_2\},$ then $$\sgn(\zeta_1;\a)\sgn(\zeta_2;\a)=-\sgn(\zeta_1;\b)\sgn(\zeta_2;\b).$$
 \smallskip

\noindent{\bf Fact 2.} If $\a,\b\in R$ such that  $\supp(\a)\neq \supp(\b)$ and  $\supp(\a)\cap \supp(\b)=\{\ep\},$ then $\a-\b\in R$ if and only if
 $\sgn(\ep;\a)=\sgn(\ep;\b).$

\begin{lem}\label{unique}
Suppose $\Pi$ is a base of $R.$ Then $\Pi$ intersects $ \{\pm2\ep_i\}_{i=1}^m+2\bbbz\d+\d$ (resp. $\{\pm2\d_p\}_{p=1}^n+2\bbbz\d$) in at most one element.
\end{lem}
\pf  Suppose to the contrary that $\a:=2s\d_i+2k\d$ and $\b:=2s'\d_j+2k'\d,$ for some $1\leq i,j\leq n,$ $s,s'\in\{\pm1\}$  and $k,k'\in\bbbz,$ are two distinct elements of $\Pi.$ We have
\begin{align*}
(s'/2)\a-(s/2)\b&=(s'/2)(2s\d_i+2k\d)-(s/2)(2s'\d_j+2k'\d)\\
&=ss'(\d_i-\d_j)+(s'k-sk')\d\in \{0\}\cup (R\setminus\{0\}\cap \frac{1}{2}(\hbox{span}_\bbbz\Pi))
\end{align*}
 which is a contradiction.
 Also if  $\a:=2s\ep_i+2k\d+s\d$ and $\b:=2s'\ep_j+2k'\d+s'\d,$ for some $1\leq i,j\leq m,$ $s,s'\in\{\pm1\}$  and $k,k'\in\bbbz,$ are two distinct elements of $\Pi,$ then as above
 \[(s'/2)\a-(s/2)\b\in \{0\}\cup (R\setminus\{0\}\cap \frac{1}{2}(\hbox{span}_\bbbz\Pi))
\]
 which is a contradiction. This completes the proof.
\qed

\begin{lem}\label{1-3}
Suppose that  $\a,\b_1,\b_2,\b_3\in R$ such that for $1\leq i\leq 3,$ $(\a,\b_i)_*\neq 0$ and with respect to $\fm_*,$ $\b_1,\b_2$ and $\b_3$ are mutually orthogonal.   Then there are $1\leq i,j\leq 3$ with $i\neq j$ such that $\supp(\b_i)=\supp(\b_j)$ and $|\supp(\b_i)|=|\supp(\b_j)|=2.$
\end{lem}
\pf
Since  for each $1\leq i\leq 3,$  $(\a,\b_i)_*\neq 0,$ we have  $\supp(\a)\cap \supp(\b_i)\neq \emptyset.$ But $|\supp(\a)|\leq 2,$ so there are $i,j$ with  $i\neq j$ such that $\supp(\a)\cap \supp(\b_i)\cap \supp(\b_j)\neq \emptyset.$ This together with Fact~1 completes the proof. \qed

\begin{cor}\label{not-E}
There is no  subset $S$ of $R$ such that $(\sspan_\bbbr S,\fm_*,S )$ is an irreducible  finite root system of type $E_{6,7,8}.$
\end{cor}
\pf  To the contrary, assume there is a subset $S$ of $R$ such that  $(\sspan_\bbbr S,\fm_*,S)$ is a finite root system of type $E_{6,7,8}.$ The  Dynkin diagram of $S$ has a sub-diagram as
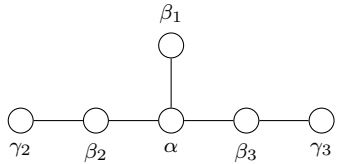
\begin{center}
\begin{tikzpicture}
\node[draw, circle, label=below:{\footnotesize $\gamma_2$}] (c1) at (9,0) {};
\node[draw, circle, label=below:{\footnotesize $\b_2$}] (c3) at (10,0) {};
\node[draw, circle, label=below:{\footnotesize $\a$}] (c4) at (11,0) {};
\node[draw, circle, label=above:{\footnotesize $\b_1$}] (c2) at (11,1) {};
\node[draw, circle, label=below:{\footnotesize $\b_3$}] (c5) at (12,0) {};
\node[draw, circle, label=below:{\footnotesize $\gamma_3$}] (c6) at (13,0) {};
‎{\draw[] (c1) edge node{} (c3);}‎
‎{\draw[] (c3) edge node{} (c4);}‎
‎{\draw[] (c4) edge node{} (c2);}‎
‎{\draw[] (c4) edge node{} (c5);}‎
‎{\draw[] (c5) edge node{} (c6);}‎
‎\end{tikzpicture}‎
\end{center}
 By Lemma \ref{1-3}, there are $1\leq i\neq j\leq 3$ with $\supp(\b_i)=\supp(\b_j)$ and  $|\supp(\b_i)|=|\supp(\b_j)|=2.$
Without loss of generality, we assume $i=2$ and assume $\supp(\b_2)=\{\ep,\eta\}.$ Since $(\b_2,\gamma_2)_*\neq 0,$ $\supp(\gamma_2)\cap\supp(\b_2)\neq \emptyset.$ If   $\supp(\gamma_2)=\supp(\b_2),$ then $(\b_2,\gamma_2)_*\neq 0$ implies that there is $r\in\{\pm1\}$ such that
\[\sgn(\ep;\b_2)=r~\sgn(\ep;\gamma_2)\andd \sgn(\eta;\b_2)=r~\sgn(\eta;\gamma_2).\] This implies that $\b_2-r\gamma_2\in \bbbz\d.$ But $\fm_*$ is nondegenerate on $\sspan_\bbbr S,$ so $\b_2=r\gamma_2$ which is a contradiction. Therefore, $|\supp(\gamma_2)\cap\supp(\b_2)|=1$ which together with the fact that $\supp(\b_2)=\supp(\b_j)$ implies that  $(\gamma_2,\b_j)_*\neq0,$ a contradiction.\qed

\begin{lem}\label{1 in 2}
Suppose that $B$ is a subset of a base of $R.$ Assume   $\a,\b\in B$ such that  $(\a,\b)_*=0$ and  $\supp(\a)=\supp(\b).$ Then we have the following:
\begin{itemize}
\item[(a)]
$|\supp(\a)|=|\supp(\b)|=2$ and  there are  $\ep,\eta\in \supp(\a)=\supp(\b)$ such that
$\sgn(\ep;\a)=-\sgn(\ep;\b)$ and $\sgn(\eta;\a)=\sgn(\eta;\b).$
\item[(b)] Keep the same notation as in part {\rm (a)},
\begin{itemize}
\item[(i)]
$\ep\not \in \supp(\gamma)$ for all $\gamma\in B\setminus\{\a,\b\};$ in particular, for $\gamma\in B\setminus\{\a,\b\},$ $(\a,\gamma)_*\neq 0$ if and only if $(\b,\gamma)_*\neq 0.$
\item[(ii)] If $\fm_*$ is nondegenerate on $\sspan_\bbbr B,$ then there is at most one element $\gamma$  of $B\setminus\{\a,\b\}$ with $\eta\in\supp( \gamma).$

\item[(iii)] If $B$ is a base of $R$ and  $\eta\in \supp(\theta)$  for some long root $\theta$ of $B,$ then for all $\gamma\in B\setminus\{\a,\b,\theta\},$ $\eta\not\in \supp(\gamma).$
\end{itemize}
\end{itemize}
\end{lem}
\pf (a)  follows from Fact 1.

(b)(i) Choose  $p,q\in\bbbz$ and  $r_1,r_2,s_1,s_2\in\{\pm 1\}$  such that $\a=r_1\ep+r_2\eta+p\d$ and $\b=s_1\ep+s_2\eta+q\d.$ Without loss of generality, we assume $$r_1=-s_1\andd r_2=s_2.$$

We claim that $\ep\not\in \supp(\gamma)$ for all $\gamma\in B\setminus\{\a,\b\}.$ Suppose to the contrary that there is $\gamma\in B\setminus\{\a,\b\}$ such that $\ep\in \supp(\gamma).$

\noindent{\bf Case 1.} $|\supp(\gamma)|=1:$ In this case, $\gamma=2r\ep+k\d$ for some $r\in\{\pm1\}$ and $k\in\bbbz.$ If $r=r_1,$ then $\a-\gamma\in R\cap (B-B)$ which is a contradiction as $B$ is a subset of  a  base of $R$ . Also if $r=-r_1=s_1,$ then $\b-\gamma\in R\cap (B-B)$ which is again a contradiction.

\noindent{\bf Case 2.} $|\supp(\gamma)|=2:$ Suppose that  $\supp(\gamma)=\{\ep,\zeta\}.$ So  $\gamma=r\ep+s\zeta+k\d$ for some $r,s\in\{\pm1\}$ and $k\in\bbbz.$ We first assume $r=r_1.$ If $\zeta\neq \eta,$ then we have $\gamma-\a=s\zeta-r_2\eta+(k-p)\d\in R\cap (B-B)$ which is a contradiction. Also if   $\zeta=\eta,$ two cases can happen: $r_2=s_2=s$ or $r_2=s_2=-s.$ In the former case, we have $(k-p)\d=\gamma-\a\in R\cap (B-B)$ which is a contradiction and in the latter case, we have   $s_1\ep+s\eta+(-p+q+k)\d=-\a+\b+\gamma$ which is an element of $R$ written as  a linear combination of the elements of a base of $R$ with coefficients of opposite sign, a  contradiction.
Finally, if $r\neq r_1,$ then $r=s_1$ and the same argument as above with changing the role of $\a$ and $\b$  gives  a contradiction. This completes the proof.

(b)(ii) Assume  $\gamma_1$ and $\gamma_2$ are two distinct elements of $B\setminus\{\a,\b\}$ with $\eta\in \supp(\gamma_1)\cap \supp(\gamma_2).$
 Since $\b-\gamma_1,\b-\gamma_2\not\in R,$ using  part  (b)(i) together with Fact 2,
we have
\begin{equation}\label{equality}
\sgn(\eta;\gamma_1)=-\sgn(\eta;\b)\andd \sgn(\eta;\gamma_2)=-\sgn(\eta;\b).
\end{equation}
If $\supp(\gamma_1)\neq \supp(\gamma_2),$ then (\ref{equality}) implies that $\gamma_1-\gamma_2\in R$ which is a contradiction.
So  $\supp(\gamma_1)=\supp(\gamma_2).$ Suppose    $\supp(\gamma_1)=\supp(\gamma_2)=\{\eta,\xi\}.$ If  $\sgn(\xi;\gamma_1)=\sgn(\xi;\gamma_2),$ (\ref{equality}) implies that $\gamma_1-\gamma_2\in \bbbz\d\sub R\cap (B-B)$ which is a contradiction. Therefore, $\sgn(\xi;\gamma_1)=-\sgn(\xi;\gamma_2),$ so we have
\[\b+\gamma_1+\gamma_2+\a\in\bbbz\d\cap \sspan_\bbbz B.\] This is a contradiction as  $\fm_*$ is nondegenerate on $\sspan_\bbbr B.$
This completes the proof.

(b)(iii)
To the contrary, assume $\eta\in\supp(\gamma)$ for some $\gamma\in B\setminus\{\a,\b,\theta\}.$  We first mention that using Lemma \ref{unique}, $\supp(\theta)\neq \supp(\gamma).$ Since  $\gamma-\theta,\theta-\b,\a-\gamma\not\in R$ and   $\ep\not\in \supp(\gamma)\cup \supp(\theta)$ (see part (b)(i)), we have
\begin{align*}
\sgn(\eta;\gamma)\stackrel{\hbox{\footnotesize Fcat 2}}{=\joinrel=}-\sgn(\eta;\a)\stackrel{\hbox{\footnotesize part (a)}}{=\joinrel=}
-\sgn(\eta;\b)\stackrel{\hbox{\footnotesize Fcat 2}}{=\joinrel=}\sgn(\eta;\theta)\stackrel{\hbox{\footnotesize Fcat 2}}{=\joinrel=}-\sgn(\eta;\gamma)
\end{align*}
which is a contradiction.
\qed

\bigskip

In what follows we shall show that up to $\pm1$-multiple,  the bases introduced in Proposition \ref{bases-last} are the only bases of $R$ of type $A(2m-1,2n-1)^{(2)},$ $(m,n)\neq (1,1).$

\begin{deft}
{\rm
Recall $\theta_i$'s and $\Pi$ from Theorem \ref{main} and assume $\supp(\theta_i)=\{\upsilon_i\}.$ We say $\Pi$ is {\it fine} if $\sgn(\upsilon_i;\theta_i)=1.$
}
\end{deft}

{\bf \emph{ From now on, we denote the bases introduced in the last four rows of Table \ref{table2} respectively by $\boldsymbol{B^1},$ $\boldsymbol{B^2},$ $\boldsymbol{B^3}$ and
$\boldsymbol{B^4}.$}}

\begin{lem}
If  $\Pi$ is  a base of the from $B^i$ ($1\leq i\leq 4$), then there is $w$ belonging to the quasi-Weyl group $W$ of $R$ such that $w(B^i)$ is a fine base of the from $B^i.$
\end{lem}
\pf
Suppose $\theta_i$'s are as in Theorem \ref{main} and assume $\supp(\theta_i)=\{\upsilon_i\}.$  For $1\leq t\leq \ell,$ set
\[r_t:=\left\{
\begin{array}{ll}
r_{2\upsilon_t}& \hbox{if  $\upsilon_t\in\{\d_i\mid 1\leq i\leq n\}$ and $\sgn(\upsilon_t;\theta_t)=-1$},\\
r_{2\upsilon_t+\d}& \hbox{if $\upsilon_t\in \{\ep_i\mid 1\leq i\leq m\}$ and  $\sgn(\upsilon_t;\theta_t)=-1$},\\
{\rm id} & \hbox{otherwise.}
\end{array}
\right.\]
Then $r_\ell\cdots r_1(R)\sub R;$ in particular, $\Pi':=r_1\cdots r_\ell(\Pi)$ is a base of $R;$ see Lemma \ref{base} and Corollary \ref{belong-to}.  Moreover, $\Pi'$ is of the form $B^i.$ Setting $\theta'_i:=r_1\cdots r_\ell(\theta_i),$ we have $\supp(\theta'_i)=\{\upsilon_i\}$ with
$\sgn(\upsilon_i;\theta'_i)=1.$  This completes the proof.
\qed

\begin{deft}\label{admissible}{\rm
Recall $\theta_i$'s and  $k_i$'s from  Theorem \ref{main} and  assume  $\Pi$ is fine of the form $B^j$  for some $1\leq j\leq 4.$ Suppose for each $1\leq i\leq \ell,$   $\supp(\theta_i)=\{\upsilon_i\}.$ In this case, we denote $\Pi$ by $(j; \upsilon_1,\ldots,\upsilon_\ell;k_1,\ldots,k_\ell).$
\begin{itemize}
\item[(i)] We say an element $w$ of the quasi-Weyl group $W$ of $R$ {\it preserves the nature} of $\Pi$ if $w(\Pi)$ is a base of $R$ of the from $B^j$
with $\supp(\theta_i)=\supp(w(\theta_i))$ for all $i.$
\item[(ii)]
We say $\Pi$ is  $t$-{\it admissible} ($1\leq t\leq \ell$) if $k_1=\cdots=k_t=0;$ in this case, we denote $\Pi$ by $(j; \upsilon_1,\ldots,\upsilon_\ell;k_{t+1},\ldots,k_\ell).$
\end{itemize}}
\end{deft}

\begin{lem}\label{basic}
\begin{itemize}
\item[(i)]
Suppose $B=(j;\upsilon_1,\ldots,\upsilon_\ell;k_{t+1},\ldots,k_\ell)$ is a $t$-admissible base for some $t\in\{1,\ldots,\ell-2\}.$ Then there is $w\in W$ such that  $w$ preserves the nature of $\Pi$   and $w(B)$ is $(t+1)$-admissible.
\item[(ii)] If  $2\leq j\leq 4$ and  $B=(j;\upsilon_1,\ldots,\upsilon_\ell;k_{2},\ldots,k_\ell)$ is a $1$-admissible base, then   there is $w\in W$
such that  $w$ preserves the nature of $\Pi$ and $w(B)$ is $\ell$-admissible.
\end{itemize}
\end{lem}
\pf
(i)  Set
\[w:=r_{\upsilon_{t+1}-\upsilon_{t+2}}r_{\upsilon_{t+1}-\upsilon_{t+2}+k_{t+1}\d},\] then by Corollary \ref{belong-to}, $w(B)$ is  a base of $R.$ Moreover,
\[w(\upsilon_{t+1})=\upsilon_{t+1}-k_{t+1}\d,\;w(\upsilon_{t+2})=\upsilon_{t+2}+k_{t+1}\d\andd w(\upsilon_i)=\upsilon_i\quad(i\neq t+1,t+2);\]  in particular,
\[w(\theta_i)=\upsilon_i\quad (1\leq i\leq t+1)\andd \supp(w(\theta_r))=\supp(\theta_r)\quad (1\leq r\leq \ell).\]This completes the proof.

(ii) It follows from (i) together with an induction process that there is   $w\in W$
such that  $w$ preserves the nature of $B$   and $B':=w(B)$ is an $(\ell-1)$-admissible base. So $B'=(j; \upsilon_1,\ldots,\upsilon_\ell;k'_{\ell}).$
We continue the proof in the following cases:

\noindent {$\bullet~\boldsymbol{j=2.}$} In this case, we have $\upsilon_\ell\in \{\d_p\mid 1\leq p\leq n\}$ and
 $$\hbox{\small$ B'=\{-2\upsilon_1,\upsilon_i-\upsilon_{i+1},\upsilon_{\ell-1}-\upsilon_{\ell}-k'_{\ell}\d,
\upsilon_{\ell-1}+\upsilon_{\ell}+(k'_\ell+1)\d\mid 1\leq i\leq \ell-2\}.$}$$
If $k'_\ell$ is odd, set
\[w_0:=r_{2\upsilon_{\ell}+(k'_\ell+1)\d},\] then using Lemma \ref{base} and Corollary \ref{belong-to}, $w_0(B')$ is a base of $R,$
\[w_0(\upsilon_\ell)=-\upsilon_\ell-k'_\ell\d\andd w_0(\upsilon_i)=\upsilon_i\quad(1\leq i\leq \ell-1);\] in particular,
\begin{align*}
w_0(\theta_{\ell-1}-\theta_\ell)&=w_0(\upsilon_{\ell-1}-\upsilon_{\ell}-k'_\ell\d)=\upsilon_{\ell-1}+\upsilon_{\ell}+\d
\\
w_0(\theta_{\ell-1}+\theta_\ell+\d)&=w_0(\upsilon_{\ell-1}+\upsilon_{\ell}+(k'_{\ell}+1)d)=\upsilon_{\ell-1}-\upsilon_{\ell}
\end{align*} and we are done in this case.  Also if $k'_{\ell}$ is even, we set
\[w_0:=r_{2\upsilon_{\ell}}r_{2\upsilon_\ell+k'_{\ell}\d}.\] We have again that $w_0(B')$ is a base of $R$ and moreover, we have
\[w_0(\upsilon_\ell)=\upsilon_\ell-k'_\ell\d\andd w_0(\upsilon_i)=\upsilon_i\quad(1\leq i\leq \ell-1);\] in particular,
\begin{align*}
w_0(\theta_{\ell-1}-\theta_\ell)&=w_0(\upsilon_{\ell-1}-\upsilon_{\ell}-k'_\ell\d)=\upsilon_{\ell-1}-\upsilon_{\ell}
\\
w_0(\theta_{\ell-1}+\theta_\ell+\d)&=w_0(\upsilon_{\ell-1}+\upsilon_{\ell}+(k'_{\ell}+1)d)=\upsilon_{\ell-1}+\upsilon_{\ell}+\d.
\end{align*}This completes the proof for $j=2.$.

\noindent {$\bullet~\boldsymbol{j=3.}$} In this case, we have $\upsilon_\ell\in \{\ep_i\mid 1\leq i\leq m\}$ and
 $$\hbox{\small$ B'=\{-2\upsilon_1+\d,\upsilon_i-\upsilon_{i+1},\upsilon_{\ell-1}-\upsilon_{\ell}-k'_{\ell}\d,
\upsilon_{\ell-1}+\upsilon_{\ell}+k'_\ell\d\mid 1\leq i\leq \ell-2\}.$}$$
Set
\[w_0:=\left\{\begin{array}{ll}
r_{2\upsilon_{\ell}+k'_\ell\d} &\hbox{if $k'_\ell$ is odd,}\\
r_{2\upsilon_{\ell}+\d}r_{2\upsilon_{\ell}+(k'_\ell+1)\d} &\hbox{if $k'_\ell$ is even,}\end{array}\right.\] then as above  $w_0(B')$ is   an $\ell$-admissible  base.

\noindent {$\bullet~\boldsymbol{j=4.}$} In this case, we have $\upsilon_1\in\{\d_p\mid 1\leq p\leq m\},\upsilon_\ell\in \{\ep_i\mid 1\leq i\leq m\}$ and
 $$\hbox{\small$ B'=\{-2\upsilon_1,\upsilon_i-\upsilon_{i+1},\upsilon_{\ell-1}-\upsilon_{\ell}-k'_{\ell}\d,
2\upsilon_{\ell}+2k'_\ell\d+\d\mid 1\leq i\leq \ell-2\}.$}$$
Setting
\[w_0:=\left\{\begin{array}{ll}
r_{2\upsilon_1}r_{\upsilon_1-\upsilon_\ell}r_{\upsilon_1+\upsilon_\ell}r_{2\upsilon_{\ell}+k'_\ell\d} &\hbox{if $k'_\ell$ is odd,}\\
r_{2\upsilon_{\ell}+\d}r_{2\upsilon_{\ell}+(k'_\ell+1)\d} &\hbox{if $k'_\ell$ is even,}\end{array}\right.\] we get that
\[w_0(\upsilon_\ell)=\upsilon_\ell-k'_\ell\d\andd w_0(\upsilon_i)=\upsilon_i\quad(1\leq i\leq \ell-1)\] and that
 $w_0(B')$ is   an $\ell$-admissible  base; see Lemma \ref{base} and Corollary \ref{belong-to}.
\qed

\begin{cor}\label{cor2}
If $B$ and $B'$ are two bases of $R$ of the form $B^j$ for some $2\leq j\leq 4,$ then $B$ and $B'$ are conjugate under the quasi-Weyl group.
\end{cor}
\pf Using Lemma \ref{basic}(ii), without loss of generality, we assume $B$ and $B'$ are $\ell$-admissible. So there are distinct elements $\upsilon_1,\ldots,\upsilon_\ell$ as well as distinct elements $\upsilon'_1,\ldots,\upsilon'_\ell$ of $\{\ep_i,\d_p\mid 1\leq i\leq m,1\leq p\leq n\}$ such that
\[B=(j;\upsilon_1,\ldots,\upsilon_\ell)\andd B'=(j;\upsilon'_1,\ldots,\upsilon'_\ell).\]

\noindent{\bf Case 1.} $\upsilon_1=\upsilon'_1$ and $\upsilon_\ell=\upsilon'_{\ell}:$
Define
\[w_2:=\left\{\begin{array}{ll}
{\rm id}& \hbox{if $\upsilon_2=\upsilon'_2$}\\
r_{\upsilon_2-\upsilon'_2}&\hbox{if $\upsilon_2\neq \upsilon'_2.$}
\end{array}\right.
\] We next assume $w_2,\ldots,w_t$ ($1\leq t\leq \ell-2$) have been defined and set
\[\a_{t+1}:=w_t\cdots w_2(\upsilon'_{t+1})\andd w_{t+1}:=\left\{\begin{array}{ll}
{\rm id}& \hbox{if $\upsilon_{t+1}=\a_{t+1}$}\\
r_{\upsilon_{t+1}-\a_{t+1}}&\hbox{if $\upsilon_{t+1}\neq \a_{t+1}.$}
\end{array}\right.
 \]
We claim that
\begin{equation}
\label{equal}
w_{t}\cdots w_2(\upsilon'_1)=\upsilon_1\andd w_{t}\cdots w_2(\upsilon'_j)=\upsilon_j\quad\quad(2\leq j\leq t\leq \ell).
\end{equation}
We use an induction process to do this; by definition, we have  \[w_2(\upsilon'_1)=w_2(\upsilon_1)=\upsilon_1\andd  w_2(\upsilon'_2)=\upsilon_2.\] Next assume $2\leq t\leq \ell-1$ and that
$w_{t}\cdots w_2(\upsilon'_j)=\upsilon_j$ for $1\leq j\leq t.$ We shall show
\[w_{t+1}\cdots w_2(\upsilon'_j)=\upsilon_j\quad (1\leq j\leq t+1).\]
Suppose $1\leq j\leq t.$ Since   $\upsilon_1,\ldots,\upsilon_\ell$ are distinct, we have
 $\upsilon_j\neq \upsilon_{t+1}.$ Also  $\upsilon_j\neq \a_{t+1}$ as otherwise
by the induction hypothesis, we have
\[w_t\cdots w_2(\upsilon'_{t+1})=\a_{t+1}=\upsilon_j=w_t\cdots w_2(\upsilon'_j)\] which in turn implies that $\upsilon'_{t+1}=\upsilon'_j,$ a contradiction. Therefore, $w_{t+1}(\upsilon_j)=\upsilon_j$ and so we have using induction hypothesis that
\begin{align*}
w_{t+1}w_t\cdots w_2(\upsilon'_j)=\left\{
\begin{array}{ll}
w_{t+1}(\upsilon_j)=\upsilon_j&1\leq j\leq t \\
w_{t+1}(\a_{t+1})=\upsilon_{t+1}=\upsilon_j & j=t+1.
\end{array}\right.
\end{align*}
 This completes the proof of (\ref{equal}). Now set $w:=w_{\ell}\cdots w_2,$ we have $w(B')=B$ as we desired.

\noindent{\bf Case 2.} General Case: We recall Lemma \ref{base} and Corollary \ref{belong-to} and set
\[\omega_1:=\left\{
\begin{array}{ll}
{\rm id}& \hbox{if $\upsilon_1=\upsilon'_1$}\\
r_{\upsilon_1-\upsilon_1'}& \hbox{if $\upsilon_1\neq \upsilon'_1.$}
\end{array}\andd \Pi:=\omega_1(B').
\right.
\] Then  for $\upsilon''_i:=\omega_1(\upsilon'_i)\in\{\ep_i,\d_p\mid 1\leq i\leq m,1\leq p\leq n\},$ we have
\[\Pi=(j;\upsilon''_1,\ldots,\upsilon''_\ell)\;\;\hbox{with} \;\;\upsilon''_1=\upsilon_1.\]
Next we set
\[\omega_2:=\left\{
\begin{array}{ll}
{\rm id}& \hbox{if $\upsilon_\ell=\upsilon''_\ell$}\\
r_{\upsilon_\ell-\upsilon''_\ell}& \hbox{if $\upsilon_\ell\neq \upsilon''_\ell$}
\end{array}\andd \Pi':=\omega_2(\Pi).
\right.
\] Then  for $\upsilon'''_i:=\omega_2(\upsilon''_i)\in\{\ep_i,\d_p\mid 1\leq i\leq m,1\leq p\leq n\},$ we have
\[\Pi'=(j;\upsilon'''_1,\ldots,\upsilon'''_\ell)\;\;\hbox{with} \;\;\upsilon'''_1=\upsilon_1\andd \upsilon'''_\ell=\upsilon_\ell .\]
So we are done using Case 1.\qed

\begin{lem}\label{types}
Suppose that  $B$ is a nonempty subset of $R^\times$ with $|B\cap R_{lg}|\leq 1;$ see (\ref{long}). Assume $W_B$ is the subgroup  of quasi-Weyl group  $W$ of $R$ generated by the quasi-reflections based on the elements of $B.$   Set $S:=W_BB$ and assume $(\sspan_\bbbr S,\fm_*,S)$ is an irreducible finite root system with base $B$, then it is  of one of   types $A_\frak{r}$ ($\frak{r}\geq1$), $D_\frak{r}$ ($\frak{r}\geq4$) or $C_\frak{r}$ ($\frak{r}\geq2$).
\end{lem}
\pf If $B\cap R_{lg}=\emptyset,$ then  $S\sub R.$   So using Corollary \ref{not-E}, $S$  is not of type $E_{6,7,8}.$
Since  bases of finite root systems of types $F_4$ and $B_\frak{r}$ ($\frak{r}\geq3$) contain at least two long roots and $B$ contains at most one long root, we get that $S$ is not neither of type $B_\frak{r}$ ($\frak{r}\geq3$) nor of type $F_4.$  Also as the  ratio of the lengths of the roots in $S$ is $1,$ $1/2$ or $2,$ $S$ is not also of type $G_2.$ \qed

\begin{lem}\label{have-delta}Suppose that $\Pi$ is a base of $R$ and $\a\in \Pi$ is a long root. Suppose $B\sub \Pi\setminus\{\a\}$ is such that for $S:=W_BB,$ $(\sspan_\bbbr S,\fm_*,S)$ is an irreducible finite root system of rank $\frak{r}$ and $B$ is a base of $S.$
\begin{itemize}
\item[(i)]  If $S$ is of type $D_\frak{r}$ $(\frak{r}\geq 4)$, then there are  distinct elements $\dot\zeta_1,\ldots,\dot\zeta_{ \frak{r}}\in \{\ep_t,\d_p\mid 1\leq t\leq m,1\leq p\leq n\}$   and $\a_1,\ldots,\a_{ \frak{r}}\in R$ with
\begin{itemize}
\item[$\bullet$] {\small $\supp(\a_i)=\{\dot\zeta_{i},\dot\zeta_{i+1}\}$ for $1\leq i\leq  \frak{r}-1,$}
\item[$\bullet$] {\small $\supp(\a_{ \frak{r}-1})=\supp(\a_ \frak{r})$ and  $(\a_{ \frak{r}-1},\a_ \frak{r})_*=0,$}
\item[$\bullet$] {\small for  $1\leq i\leq  \frak{r} -1 ,$    $\sgn(\dot\zeta_{i+1};\a_i)=-\sgn(\dot\zeta_{i+1};\a_{i+1}),$  }
\item[$\bullet$] {\small $\sgn(\dot\zeta_{\frak{r}-1};\a_{\frak{r}-1})=\sgn(\dot\zeta_{\frak{r}-1};\a_{\frak{r}}),$}
\item[$\bullet$] {\small $B=\{\a_1,\ldots,\a_{ \frak{r}}\}.$}
\end{itemize}

\item[(ii)]
If $S$ is of type  $C_{\frak r}$ $(\frak{r}\geq 2),$ then there are  distinct elements $\dot\zeta_1,\ldots,\dot\zeta_{ \frak{r}}\in \{\ep_t,\d_p\mid 1\leq t\leq m,1\leq p\leq n\}$  and $\a_1,\ldots,\a_{ \frak{r}}\in R$ with
\begin{itemize}
\item[$\bullet$] {\small $\supp(\a_i)=\{\dot\zeta_{i},\dot\zeta_{i+1}\}$ for $1\leq i\leq  \frak{r}-1,$}
\item[$\bullet$] {\small  $\supp(\a_{ \frak{r}})=\{\dot\zeta_{ \frak{r}}\},$}
\item[$\bullet$] {\small for  $1\leq i\leq  \frak{r}  ,$    $\sgn(\dot\zeta_{i+1};\a_i)=-\sgn(\dot\zeta_{i+1};\a_{i+1}),$}
\item[$\bullet$] {\small $B=\{\a_1,\ldots,\a_{ \frak{r}}\}.$}
\end{itemize}
\item[(iii)]
If $S$ is of type  $A_{\frak r}$ for some   $\frak{r}\neq 1,3$ or $S\cap R_{lg}=\emptyset$ and $\frak{r}=1,$ then there are  distinct elements $\dot\zeta_1,\ldots,\dot\zeta_{ \frak{r}+1}\in \{\ep_t,\d_p\mid 1\leq t\leq m,1\leq p\leq n\}$  and $\{\a_1,\ldots,\a_{ \frak{r}}\}$ with
\begin{itemize}
\item[$\bullet$] {\small $\supp(\a_i)=\{\dot\zeta_{i},\dot\zeta_{i+1}\}$ for $1\leq i\leq  \frak{r},$}
\item[$\bullet$] {\small for  $1\leq i\leq  \frak{r}  ,$    $\sgn(\dot\zeta_{i+1};\a_i)=-\sgn(\dot\zeta_{i+1};\a_{i+1}),$}
\item[$\bullet$] {\small $B=\{\a_1,\ldots,\a_{ \frak{r}}\}.$}
\end{itemize}
\item[(iv)]
If $S$ is of type  $A_{3},$ then one of the following happens:
\begin{itemize}
\item[(a)] There are $\a_1,\a_2,\a_3\in R$  and  distinct elements $\dot\zeta_1,\dot\zeta_2,\dot\zeta_3,\dot\zeta_4\in \{\ep_t,\d_p\mid 1\leq t\leq m,1\leq p\leq n\}$   with
\begin{itemize}
\item[$\bullet$] {\small $\supp(\a_i)=\{\dot\zeta_{i},\dot\zeta_{i+1}\}$ for $1\leq i\leq  3,$}
\item[$\bullet$] {\small for  $1\leq i\leq  3  ,$    $\sgn(\dot\zeta_{i+1};\a_i)=-\sgn(\dot\zeta_{i+1};\a_{i+1}),$}
\item[$\bullet$] {\small $B=\{\a_1,\a_2,\a_3\}.$}
\end{itemize}
\item[(b)] There are  distinct elements $\dot\zeta_1,\dot\zeta_2,\dot\zeta_3\in \{\ep_t,\d_p\mid 1\leq t\leq m,1\leq p\leq n\}$  and $\a_1,\a_2,\a_3\in R$ with
\begin{itemize}
\item[$\bullet$] {\small $\supp(\a_i)=\{\dot\zeta_{i},\dot\zeta_{i+1}\}$ for $1\leq i\leq  2,$}
\item[$\bullet$] {\small $\supp(\a_2)=\supp(\a_3)$ and  $(\a_{ 2},\a_ 3)_*=0,$}
\item[$\bullet$] {\small for  $1\leq i\leq  2,$    $\sgn(\dot\zeta_{i+1};\a_i)=-\sgn(\dot\zeta_{i+1};\a_{i+1}),$  }
\item[$\bullet$] {\small $\sgn(\dot\zeta_{2};\a_{2})=\sgn(\dot\zeta_{2};\a_{3}),$}
\item[$\bullet$] {\small $B=\{\a_1,\a_2,\a_{3}\}.$}
\end{itemize}
\end{itemize}
\end{itemize}
\end{lem}
\pf
(i) By our assumption,  $B$ is a base of a finite root system of type $D_ \frak{r}.$ Suppose  $B=\{\a_1,\ldots,\a_{ \frak{r}}\}$ with  the corresponding Dynkin diagram
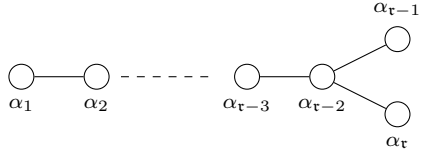
\begin{center}
\begin{tikzpicture}
\node[draw, circle, label=below:{\footnotesize $\a_1$}] (c1) at (9,0) {};
\node[draw, circle, label=below:{\footnotesize $\a_2$}] (c2) at (10,0) {};
\node (ce1) at (10.2,0) {};%empty node
\node (ce2) at (11.7,0) {};%empty node
  { \draw[dashed](ce1) -- (ce2);}
\node[draw, circle, label=below:{\footnotesize $\a_{ \frak{r}-3}$}] (c3) at (12,0) {};
\node[draw, circle, label=below:{\footnotesize $\a_{ \frak{r}-2}$}] (c4) at (13,0) {};
\node[draw, circle, label=above:{\footnotesize $\a_{ \frak{r}-1}$ }] (c5) at (14,.5) {};
\node[draw, circle, label=below:{\footnotesize $\a_{ \frak{r}}$}] (c6) at (14,-.5) {};
{\draw[] (c1) edge node{} (c2);}
{\draw[] (c3) edge node{} (c4);}
{\draw[] (c4) edge node{} (c5);}
{\draw[] (c4) edge node{} (c6);}
\end{tikzpicture}
\end{center}
and note that by Lemma \ref{unique}, for each $i,$ $|\supp(\a_i)|=2.$

\begin{center}{\bf Convention.} If $ \frak{r}=4,$ changing the indices if necessary and using Lemma \ref{1-3}, we assume $\supp(\a_{ \frak{r}})=\supp(\a_{ \frak{r}-1}).$
\end{center}

\noindent{\bf Step 1.} We have $\supp(\a_\frak{r})=\supp(\a_{ \frak{r}-1}):$ It follows from  the above convention, Lemmas \ref{1-3} and \ref{1 in 2}(b)(i) together with the fact that if $\frak{r}>4,$ we have  $(\a_{\frak{r}-4},\a_{\frak{r}})_*,(\a_{\frak{r}-4},\a_{\frak{r}-1})_*=0 $ while $(\a_{\frak{r}-4},\a_{\frak{r}-3})_*\neq 0.$

\medskip

\noindent{\bf Step 2.} For $1\leq i<j\leq \frak{r}-2$, we have  $\supp(\a_i)\cap \supp(\a_{j+1})=\emptyset:$ Suppose to the contrary that
$1\leq i<j\leq \frak{r}-2$ and  $\supp(\a_i)\cap \supp(\a_{j+1})\neq \emptyset.$ Since  $(\a_i,\a_{j+1})_*=0,$  we have $\supp(\a_i)= \supp(\a_{j+1})$ and using   Lemma \ref{1 in 2}(b)(i),  for $\gamma\in B\setminus\{\a_i,\a_{j+1}\},$ $(\gamma,\a_i)_*\neq 0$ if and only if $(\gamma,\a_{j+1})_*\neq 0.$
But if $j\neq \frak{r}-2,$   $(\a_i,\a_{j+2})_*= 0$ while  $(\a_{j+1},\a_{j+2})_*\neq  0,$ a contradiction. Also if $j=\frak{r}-2,$ using Step 1,  we have
$\supp(\a_i)=\supp(\a_{j+1})=\supp(\a_{j+2})$ which contradicts Lemma \ref{1 in 2}(b)(i).

\medskip

\noindent {\bf Step 3.}
For  $1\leq i\leq  \frak{r}-2,$ $|\supp(\a_i)\cap \supp(\a_{i+1})|=1:$
 Since $(\a_i,\a_{i+1})_*\neq 0,$ we have $\supp(\a_i)\cap \supp(\a_{i+1})\neq \emptyset.$ If $|\supp(\a_i)\cap \supp(\a_{i+1})|\neq 1,$ we have
 $|\supp(\a_i)\cap \supp(\a_{i+1})|=2.$ Suppose $\supp(\a_i)= \supp(\a_{i+1})=\{\zeta,\eta\},$ then since $(\a_i,\a_{i+1})_*\neq 0,$ there is $r\in\{\pm 1\}$ with $\sgn(\zeta;\a_{i})=r~\sgn(\zeta;\a_{i+1})$ and $\sgn(\eta;\a_{i})=r~ \sgn(\eta;\a_{i+1}).$ But this implies that $\a_i-r\a_{i+1}\in\bbbz\d.$ Since  $\fm_*$ is nondegenerate on $\sspan_\bbbr B,$ we have  $\a_i-r\a_{i+1}=0$ which is a contradiction.

  \medskip

\noindent{\bf Step 4.} Use Steps 1--3 to   pick  distinct elements  $\dot\zeta_{i}\in\{\ep_t,\d_p\mid 1\leq t\leq m,1\leq j\leq n\}$ ($1\leq i\leq \frak{r}$) with
$\dot\zeta_i\in\supp(\a_i)$  and for $1\leq i\leq \frak{r}-2,$  $\supp(\a_i)\cap \supp(\a_{i+1})=\{\dot\zeta_{i+1}\}.$ For $1\leq i\leq \frak{r}-2,$  $\sgn(\dot\zeta_{i+1};\a_{i})=-\sgn(\dot\zeta_{i+1};\a_{i+1})$ as otherwise, $\a_i-\a_{i+1}\in R$ which is a contradiction (see Fact 2). Similarly, $\sgn(\dot\zeta_{\frak{r}-1};\a_{\frak{r}-2})=-\sgn(\dot\zeta_{\frak{r}-1};\a_{\frak{r}}).$ This implies that  $\sgn(\dot\zeta_{\frak{r}-1};\a_{\frak{r}-1})=\sgn(\dot\zeta_{\frak{r}-1};\a_{\frak{r}})$ and so $\sgn(\dot\zeta_{\frak{r}};\a_{\frak{r}-1})=-\sgn(\dot\zeta_{\frak{r}};\a_{\frak{r}})$ (see Fact 1).

(ii) is similarly proved.

(iii),(iv) Suppose  $B=\{\a_1,\ldots,\a_{\frak{r}}\}$ such that  the corresponding Dynkin diagram of $S$ is
as follows:
 
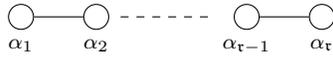
\begin{center}
\begin{tikzpicture}
\node[draw, circle, label=below:{\footnotesize $\a_1$}] (c1) at (9,0) {};
\node[draw, circle, label=below:{\footnotesize $\a_2$}] (c2) at (10,0) {};
\node (ce1) at (10.2,0) {};%empty node
\node (ce2) at (11.7,0) {};%empty node
  { \draw[dashed](ce1) -- (ce2);}
\node[draw, circle, label=below:{\footnotesize $\a_{\frak{r}-1}$}] (c3) at (12,0) {};
\node[draw, circle, label=below:{\footnotesize $\a_{\frak{r}}$}] (c4) at (13,0) {};
{\draw[] (c1) edge node{} (c2);}
{\draw[] (c3) edge node{} (c4);}
\end{tikzpicture}
\end{center}
As in Step 3 of the proof of part (i), for $1\leq i\leq \frak{r}+1,$  there is $\dot\zeta_i\in\{\ep_j,\d_p\mid 1\leq j\leq m,1\leq p\leq n\}$ with
$\dot\zeta_1\in\supp(\a_1),$ $\dot\zeta_{\frak{r}+1}\in\supp(\a_{\frak{r}})$ and for $1\leq i\leq \frak{r}-1,$ $\supp(\a_i)\cap  \supp(\a_{i+1})=\{\dot\zeta_{i+1}\}$  and as in the proof of Step 4 of part (i), $\sgn(\dot\zeta_{i+1};\a_i)=-\sgn(\dot\zeta_{i+1};\a_{i+1})$ for $1\leq i\leq \frak{r}-1.$
This completes the proof if $\dot\zeta_i$'s are distinct.

Now assume $\dot\zeta_i$'s are not distinct. Then  there are  $1\leq i<i+1<j\leq \frak{r}$ with $\supp(\a_i)\cap \supp(\a_j)\neq\emptyset.$ Since $(\a_i,\a_j)_*=0,$ we have $\supp(\a_i)=\supp(\a_j).$   Assume $\supp(\a_i)=\supp(\a_j)=\{\zeta,\eta\}.$ Using Lemma \ref{1 in 2}(a),(b)(i), we assume $\sgn(\zeta;\a_i)=-\sgn(\zeta;\a_j),$ $\sgn(\eta;\a_i)=\sgn(\eta;\a_j)$  and
\begin{equation}\label{not in}
\zeta\not\in \supp(\gamma)\quad\quad(\gamma\in B\setminus\{\a_i,\a_j\}).
\end{equation}
 We know $(\a_i,\a_{i+1})_*\neq 0$ and $(\a_{j-1},\a_j)_*\neq 0,$ so
 \begin{equation*}\label{ij}
 \eta\in \supp(\a_{i+1})\cap \supp(\a_{j-1}).
 \end{equation*}
This together with  Lemma \ref{1 in 2}(b)(ii)  implies that   \[i+1=j-1.\]
Now if $i\neq 1,$ since $(\a_{i-1},\a_i)_*\neq 0,$ (\ref{not in}) implies that $\eta\in \supp(\a_{i-1})$ which is a contradiction due to Lemma \ref{1 in 2}(b)(ii) and similarly,  if $j\neq \frak{r},$  we have  $(\a_{j+1},\a_j)_*\neq 0$ which is a  contradiction as above. So $i=1$ and $j=\frak{r};$ i.e., $\frak{r}=3$ and we are done.\qed

\begin{cor}\label{cor1} Suppose that $\Pi$ is a base of $R$ and $\a\in \Pi$ is a long root. Suppose $B\sub \Pi\setminus\{\a\}$  such that for $S:=W_BB,$ $(\sspan_\bbbr S,\fm_*,S)$ is an irreducible finite root system of rank $\frak{r},$ $B$ is a base of $S$ and  $(\a,B)_*\neq \{0\}.$
\begin{itemize}
\item[(i)]  If $S$ is of type $D_\frak{r}$ $(\frak{r}\geq 4)$, then there are   $$r\in\{0,1\},\;\;k^*\in\bbbz\andd\zeta_1,\ldots,\zeta_{ \frak{r}}\in \{\pm\ep_t,\pm\d_p\mid 1\leq t\leq n,1\leq p\leq m\}+\bbbz\d$$  with   $\supp(\zeta_i)=\{\dot\zeta_i\}$  such that  $\supp(\a)=\{\dot \zeta_1\},$  $\dot\zeta_i$'s are distinct,
    $$\a=-2\zeta_1+r\d\andd B=\{\zeta_i-\zeta_{i+1},\zeta_{\frak{r}-1}+\zeta_{\frak{r}}+k^*\d\mid 1\leq i\leq \frak{r}-1\}.$$
\item[(ii)]
If $S$ is of type  $C_{\frak r}$ $(\frak{r}\geq 2),$ then there are   $$r\in\{0,1\},\;\;k^*\in\bbbz\andd\zeta_1,\ldots,\zeta_{ \frak{r}}\in \{\pm\ep_t,\pm\d_p\mid 1\leq t\leq n,1\leq p\leq m\}+\bbbz\d$$   with   $\supp(\zeta_i)=\{\dot\zeta_i\}$  such that  $\supp(\a)=\{\dot \zeta_1\},$  $\dot\zeta_i$'s are distinct,
    $$\a=-2\zeta_1+r\d\andd B=\{\zeta_i-\zeta_{i+1},2\zeta_{\frak{r}}+k^*\d\mid 1\leq i\leq \frak{r}\}.$$
 \item[(iii)]
If $S$ is of type  $A_{\frak r}$ $(\frak{r}\neq 3),$ then there are   $$r\in\{0,1\},\;\;\zeta_1,\ldots,\zeta_{ \frak{r}+1}\in \{\pm\ep_t,\pm\d_p\mid 1\leq t\leq n,1\leq p\leq m\}+\bbbz\d$$  with   $\supp(\zeta_i)=\{\dot\zeta_i\}$  such that  $\supp(\a)=\{\dot \zeta_1\},$  $\dot\zeta_i$'s are distinct,
    $$\a=-2\zeta_1+r\d\andd B=\{\zeta_i-\zeta_{i+1}\mid 1\leq i\leq \frak{r}\}.$$
\item[(iv)]
If $S$ is of type  $A_{3},$ then one of the following happens:
\begin{itemize}
\item[(a)]  There are   $$r\in\{0,1\},\;\;k^*\in\bbbz\andd\zeta_1,\zeta_2,\zeta_{3}\in \{\pm\ep_t,\pm\d_p\mid 1\leq t\leq n,1\leq p\leq m\}+\bbbz\d$$ with   $\supp(\zeta_i)=\{\dot\zeta_i\}$  such that  $\supp(\a)=\{\dot \zeta_1\},$  $\dot\zeta_i$'s are distinct,
    $$\a=-2\zeta_1+r\d\andd B=\{\zeta_1-\zeta_{2},\zeta_2-\zeta_3,\zeta_{2}+\zeta_{3}+k^*\d\}.$$

\item[(b)]  There are   $$r\in\{0,1\},\;\;\zeta_1,\zeta_2,\zeta_3,\zeta_{ 4}\in \{\pm\ep_t,\pm\d_p\mid 1\leq t\leq n,1\leq p\leq m\}+\bbbz\d$$  with   $\supp(\zeta_i)=\{\dot\zeta_i\}$  such that  $\supp(\a)=\{\dot \zeta_1\},$  $\dot\zeta_i$'s are distinct,
    $$\a=-2\zeta_1+r\d\andd B=\{\zeta_i-\zeta_{i+1}\mid 1\leq i\leq 3\}.$$
\end{itemize}

\end{itemize}
\end{cor}

\pf  (i) By Lemma \ref{have-delta}(i), there are distinct elements  $\dot\zeta_1,\ldots,\dot\zeta_\frak{r}\in\{\ep_i,\d_p\mid 1\leq i\leq m,1\leq p\leq n\},$  $k_1,\ldots,k_\frak{r}\in\bbbz$ and $\a_1,\ldots,\a_{\frak{r}}\in R$ such that $\supp(\a_i)=\{\dot\zeta_i,\dot\zeta_{i+1}\}$  ($1\leq i \leq \frak{r}-1$) and $\supp(\a_{\frak{r}})=\{\dot\zeta_{\frak{r}-1},\dot\zeta_{\frak{r}}\}.$ Also for
\[\eta_i:=\sgn(\dot\zeta_i;\a_i)\dot\zeta_i,\] $$B=\{\a_{i}:=\eta_i-\eta_{i+1}+k_{i}\d, \a_\frak{r}:=\eta_{\frak{r}-1}+\eta_{\frak{r}}+k_{\frak{r}}\d\mid 1\leq i\leq \frak{r}-1\}.$$
Since   $(\a,B)_*\neq \{0\},$ $\supp(\a)=\{\dot\zeta_j\}$ for some $1\leq j\leq \frak{r}.$ Suppose  $j\neq 1,$  then since $\a-\a_{j-1},\a-\a_j\not\in R,$ using Fact 2, we have $\sgn(\dot\zeta_j;\a)=-\sgn(\dot\zeta_j;\a_j)$ and $\sgn(\dot\zeta_j;\a)=-\sgn(\dot\zeta_j;\a_{j-1})$   which is  a contradiction. So  $\supp(\a)=\{\dot\zeta_1\}$
and as $\a-\a_1\not\in R,$ $\sgn(\dot\zeta_1;\a)=-\sgn(\dot\zeta_1;\a_1).$ So $\a=-2\eta_1+2k_0\d+r\d,$ for some $k_0\in\bbbz$ and $r\in\{0,1\}.$
 Set
\[k^*:=k_{\frak{r}}+k_{\frak{r}-1}+2\sum_{j=0}^{\frak{r}-2}k_j\andd \zeta_i:=\eta_i-\sum_{j=0}^{i-1}k_j\d\quad(1\leq i\leq \frak{r}).\] Then
we have $B=\{\a_{i}=\zeta_i-\zeta_{i+1}, \a_{\frak{r}}=\zeta_{\frak{r}-1}+\zeta_{\frak{r}}+k^*\d\mid 1\leq i\leq \frak{r}-1\}.$

(ii)-(iv) are similarly proved.\qed

%\begin{lem}\label{nece}
%Suppose that $r=0,1$ $\eta_1,\ldots,\eta_{\ell}\in\{\pm\ep_i,\pm\d_p\mid 1\leq i\leq m,1\leq p\leq n\}$ with $\eta_i\neq \pm\eta_j$ if $i\neq j$ and  $k_0,\ldots,k_{\ell}\in\bbbz$ are such that $$\hbox{\small $\Pi=\{-2\eta_1+2k_0\d+r\d,\a_{i+1}:=\eta_i-\eta_{i+1}+k_{i}\d, \a_{\ell+1}:=\eta_{\ell-1}+\eta_{\ell}+k_{\ell}\d\mid 1\leq i\leq \ell-1\}$}$$ is a base of $R.$ Then if $r=0$ (resp. $r=1$), $\eta_1,\eta_\ell\in\{\pm\d_p\mid 1\leq p\leq n\}$ (resp. $\eta_1,\eta_\ell\in\{\pm\ep_i\mid 1\leq i\leq m\}$).
%\end{lem}
%\pf Setting
%\[k^*:=k_{\ell}+k_{\ell-1}+2\sum_{j=0}^{\ell-2}k_j\andd \zeta_i:=\eta_i-\sum_{j=0}^{i-1}k_j\d\quad(1\leq i\leq \ell),\] we have
%$$\hbox{\small $\Pi=\{-2\zeta_1+r\d,\a_{i+1}:=\zeta_i-\zeta_{i+1}, \a_{\ell+1}:=\zeta_{\ell-1}+\zeta_{\ell}+k^*\d\mid 1\leq i\leq \ell-1\}.$}$$
%Since $\d\in\sspan_{\bbbz^{\geq0}}\Pi\cup \sspan_{\bbbz^{\leq0}}\Pi,$ we get that  $r+k^*=\pm1.$ This together with
%\[-2\zeta_1+r\d\in R\andd 2\zeta_\ell+(\pm1-r)\d=\a_{\ell+1}-\a_\ell\not\in R\] completes the proof.\qed
\begin{Thm}\label{finite}
Suppose that $\Pi$ is a base of $R$ and  $\a\in \Pi$ is a long root. Set  $B:=\Pi\setminus\{\a\}$ and  assume $W_B$ is the subgroup of the quasi-Weyl group  $W$ of $R$ generated by the quasi-reflections  based on the elements of $B.$ Then  $(\a,B)_*\neq \{0\}$ and for $S:=W_BB,$  $(\sspan_\bbbr S, \fm_*,S)$ is an irreducible  finite root system with base $B$,  either of type  $D_\ell$ ($\ell\geq 3$) or $C_\ell$ ($\ell\geq 2$).
\end{Thm}
\pf
We carry out the proof through the following steps:
%We know  $\a=2s\ep_i+2k^*\d+\d$ for some  $s\in\{\pm1\},$ $k^*\in\bbbz$ and $1\leq i\leq m$ or $\a=2s\d_p+2k^*\d$  for some $k^*\in\bbbz$ and $1\leq p\leq n.$ Set $B:=\Pi\setminus\{\a\}$ and   $S:=W_BB.$

\noindent{\bf Step 1.} $\bbbz\d\cap \sspan_\bbbr S=\bbbz\d\cap \sspan_\bbbr B=\{0\}:$ To  the contrary, assume $\bbbz\d\cap \sspan_\bbbr B$ contains a  nonzero element.  Then there is a nonzero integer $k$ such that $k\d\in \sspan_\bbbr B.$ But   $$2k\d\in R\sub \sspan_{\bbbz^{\geq 0}}(B\cup\{\a\})\cup \sspan_{\bbbz^{\leq 0}}(B\cup\{\a\}).$$ Since $k\d\in\sspan_\bbbr B,$ $2k\d$ is a nonzero element of $R$ belonging to $\sspan_{\bbbz^{\geq 0}}B\cup\sspan_{\bbbz^{\leq 0}}B.$ Therefore, we have one of the following:
\begin{itemize}
\item $0\neq 2k\d\in \sspan_{\bbbz^{\geq 0}}B:$ In this case, $\a-2k\d\in R$ while it is an element of $\sspan_{\bbbz^{> 0}}\{\a\}\cup \sspan_{\bbbz^{\leq 0}}(\Pi\setminus\{\a\})$ which is a contradiction.
\item $0\neq 2k\d\in \sspan_{\bbbz^{\leq 0}}B:$ In this case, $\a+2k\d\in R$ while it is an element of $\sspan_{\bbbz^{> 0}}\{\a\}\cup \sspan_{\bbbz^{\leq 0}}(\Pi\setminus\{\a\});$ again it is a contradiction.
\end{itemize}

\noindent{\bf Step 2.} If  $\a,\b\in S$ with $(\a-\b,B)_*=\{0\},$ then $\a=\b:$  To the contrary, assume there are  $\a,\b\in S$  such that $\a\neq \b$ and  $(\a-\b,B)_*=\{0\}.$  Since $\sspan_\bbbr B=\sspan_\bbbr{S},$
\[
\hbox{$(\a-\b,B)_*=\{0\}$ if and only if $(\a-\b,S)_*=\{0\}.$}
\]

\noindent{\bf Case 1.} $\a,\b\in\{\pm2\ep_i,\pm2\d_p\mid 1\leq i\leq m,1\leq p\leq n\}+\bbbz\d$:  In this case, there are $p,q\in\bbbz,$ elements  $\zeta_1,\zeta_2$ of $\{\ep_i, \d_p\mid 1\leq i\leq m,\;1\leq p\leq n\}$ and $r_1,r_2\in\{\pm1\}$
such that $\a=2r_1\zeta_1+p\d$ and $\b=2r_2\zeta_2+q\d.$  Since $\a\in S,$ $(\a-\b,\a)_*=0.$ This in turn implies that $\zeta_1=\zeta_2$ and $r_1-r_2=0,$ i.e., $(p-q)\d=\a-\b\in \sspan_\bbbr S=\sspan_\bbbr B.$ Using Step 1, we get that $p=q$. It means that $\a=\b.$

\noindent{\bf Case 2.} $\a\in\{\pm\ep_i\pm\ep_j,\pm\d_p\pm\d_q,\pm\ep_i\pm\d_p\mid 1\leq i\neq j\leq m,1\leq p\neq q\leq n\}+\bbbz\d$ and $\b\in\{2\ep_i,2\d_p\mid 1\leq i\leq m,1\leq p\leq n\}+\bbbz\d$:   In this case, there are $p,q\in\bbbz,$ $\zeta_1,\zeta_2,\zeta_3\in\{\ep_i, \d_p\mid 1\leq i\leq m,\;1\leq p\leq n\},$ with $\zeta_1\neq \zeta_2,$ and $r_1,r_2,r_3\in\{\pm1\}$ such that $\a=r_1\zeta_1+r_2\zeta_2+p\d$ and $\b=2r_3\zeta_3+q\d.$ If $\zeta_3\not\in\{\zeta_1,\zeta_2\},$
then $0=(\a-\b,\b)_*=-4$ which is a contradiction. So $\zeta_3\in\{\zeta_1,\zeta_2\},$ say e.g., $\zeta_3=\zeta_1.$ So $$0=(\a-\b,\b)_*=((r_1-2r_3)\zeta_3+r_2\zeta_2+(p-q)\d,2r_3\zeta_3+q\d)_*=2r_3(r_1-2r_3)\neq0,$$ a contradiction.

\noindent{\bf Case 3.} $\a,\b\in\{\pm\ep_i\pm\ep_j,\pm\d_p\pm\d_q,\pm\ep_i\pm\d_p\mid 1\leq i\neq j\leq m,1\leq p\neq q\leq n\}+\bbbz\d:$
Suppose that $$\a=r_1\zeta_1+r_2\zeta_2+p\d\andd \b=-r_3\zeta_3-r_4\zeta_4+q\d,$$ for some $p,q\in\bbbz,$  elements $\zeta_1,\ldots,\zeta_4$ of $\{\ep_i,\d_j\mid 1\leq i\leq m,\;1\leq p\leq n\}$ with $\zeta_1\neq \zeta_2$ and  $\zeta_3\neq \zeta_4$ and $r_1,\ldots,r_4\in\{\pm1\}.$ We first assume $\{\zeta_1,\zeta_2\}\cap\{\zeta_3,\zeta_4\}=\emptyset.$
Since $\zeta_1\in \supp(\a)$ and $\a\in S=W_BB,$ we pick $\gamma\in B$ with $\zeta_1\in \supp(\gamma).$ Since $(\a-\b,\gamma)_*=0,$ there are $r,s\in\{\pm1\},$ $k\in\bbbz$ and  $2\leq i\leq 4$ with $rs=-r_1r_i$ and $\gamma=r\zeta_1+s\zeta_i+k\d.$
Since $$\eta:=r_2\zeta_2+r_i\zeta_i+(p-rr_1k)\d=r_{\zeta_1+rs\zeta_i+rk\d}(r_1\zeta_1+r_2\zeta_2+p\d)\in W_BS\sub S,$$  we have $2=(\eta,\a-\b)_*\in(S,\a-\b)_*=\{0\},$ a contradiction. So $\{\zeta_1,\zeta_2\}\cap\{\zeta_3,\zeta_4\}\neq\emptyset,$ say $\zeta_1=\zeta_3.$

 If $r_1+r_3\neq 0,$ then  $r_1=r_3.$ So $\a-\b=2r_1\zeta_1+r_2\zeta_2+r_4\zeta_4+(p-q)\d.$ As above, we find $\gamma \in B$ with $\zeta_1\in\supp(\gamma),$ so we have $(\gamma,\a-\b)\neq 0$ which contradicts $(B,\a-\b)_*=\{0\}.$
So $r_1+r_3=0.$
 Therefore, we have  $$\a-\b=r_2\zeta_2+r_4\zeta_4+(p-q)\d.$$

 Since $\zeta_2\in \supp(\a)$ and $\a\in S=W_BB,$ we pick  $\gamma\in B$ with $\zeta_2\in \supp(\gamma).$ One knows  $(\gamma,\a-\b)=\{0\},$ so if $\zeta_2\neq \zeta_4,$
there are $k\in\bbbz$ and  $r,s\in\{\pm1\}$ with $rs=-r_2r_4$ and $\gamma=r\zeta_2+s\zeta_4+k\d.$
But $$\eta:=r_1\zeta_1+r_4\zeta_4+(p-r_2rk)\d=r_{\zeta_2+rs\zeta_4+rk\d}(r_1\zeta_1+r_2\zeta_2+p\d)\in S$$ and  $0=(\eta,\a-\b)_*=1,$ a contradiction, i.e.,  $\zeta_2=\zeta_4$ and so $$\a-\b=(r_2+r_4)\zeta_2+(p-q)\d.$$ But $\zeta_2$ belongs to the support of an element $\theta\in B.$ Since  $(\a-\b,\theta)_*=0,$ we get that $r_2+r_4=0.$ So using Step 1, we get that $p-q=0$  and so $\a=\b.$

\medskip

\noindent{\bf Step 3.} $S$ is a  finite root system: Suppose   $B=\{\b_1,\ldots,\b_\ell\}$ and define
\begin{align*}
\varphi:& S\longrightarrow \bbbz^\ell\\
&\b\mapsto(\frac{2(\b,\b_1)_*}{(\b_1,\b_1)_*},\ldots,\frac{2(\b,\a_\ell)_*}{(\b_\ell,\b_\ell)_*})
\end{align*}
We get using Step 2 that the map $\varphi$ is one to one. Moreover,  for $1\leq i\leq \ell,$ $\frac{2(\b,\b_i)_*}{(\b_i,\b_i)_*}\in\{0,\pm1,\pm2\}.$ Therefore,  $S$ is finite; in particular, $(\sspan_\bbbr S,\fm_*,S)$ is a finite root system.

\medskip

{\footnotesize \bf  Since $\boldsymbol S$ is a finite root system,  it is  a direct sum of irreducible finite root systems; say   $\boldsymbol{S=S_1\op\cdots\op S_k}.$ Setting $\boldsymbol{B_i:=B\cap S_i}$ $\boldsymbol{(1\leq i\leq k)},$ we have $\boldsymbol{S_i=W_{B_i}B_i}.$}

\medskip

%\noindent{\bf Step 4.} If $1\leq i\neq j\leq k$ and $\supp(S_i)\cap\supp(S_j)\cap  \supp(\a)\neq \emptyset,$ then $m+n\geq 3:$ Suppose to the contrary that $m+n<3,$ then $m=n=1.$ Since $|\Pi|=m+n+1=1+\sum_{i=1}^k|B_i|,$ we get $k=2\andd |B_1|=|B_2|=1.$  Assume $B_1=\{\b\}$ and $B_2=\{\gamma\}$ and $\supp(\a)=\{\eta_1\}.$ Since
%$(\a,\b)_*,(\a,\gamma)_*\neq 0,$ we have $\eta_1\in\supp(\b)\cap \supp(\gamma)$ and as  $(\b,\gamma)_*=0,$ by Lemma \ref{1 in 2}, there is   $\eta_2\in\{\ep_1,\d_1\}\setminus\{\eta_1\}$ such that $\supp(\b)=\supp(\gamma)=\{\eta_1,\eta_2\}$
%and \[\sgn(\eta_2;\b)=-\sgn(\eta_2;\gamma)\andd -\sgn(\eta_1;\a)=\sgn(\eta_1;\b)=\sgn(\eta_1;\gamma).\]
%Set \[\zeta_1:=\sgn(\eta_1;\b)\eta_1\andd \zeta_2:=\sgn(\eta_2;\b)\eta_2.\] There are $r\in\{0,1\}$ and  $t,p,q\in\bbbz$ with \[\Pi=\{\a\}\cup B_1\cup B_2=\{\a,\b,\gamma\}=\{-2\zeta_1+2t\d+r\d,\zeta_1+\zeta_2+p\d,\zeta_1-\zeta_2+q\d\}.\] Setting
%\[\dot\zeta_1:=\zeta_1-t\d,\;\;\dot\zeta_2:=\zeta_2+(t+p)\d\andd k^*:=p+q+2t,\] we have \[\Pi=\{\a=-2\dot\zeta_1+r\d,\b:=\dot\zeta_1+\dot\zeta_2,\gamma:=\dot\zeta_1-\dot\zeta_2+k^*\d\}.\] Since $\d\in\sspan_\bbbz\Pi,$ we get $r+k^*=\pm1.$ If $r=0$ (i.e., $\eta_1=\d_1$), we have $k^*=\pm1.$ On the other hand,  $2\dot\zeta_2\mp\d=2\dot\zeta_2-k^*\d=\b-\gamma\not\in R,$ we get $\eta_2=\d_1,$ a contradiction. Similarly, we get a contradiction if $r=1.$
%
%\smallskip

\noindent{\bf Step 4.} There is
a unique
$i\in\{1,\ldots,k\}$ with $\supp(B_i)\cap \supp(\a)\neq \emptyset:$
Suppose $\supp(\a)=\{\ep\}$ and $\a=2\sgn(\ep;\a)\ep+k\d$ for some  $k\in\bbbz.$ If $\ep\not\in\supp(B),$ then for each $\eta\in \supp(B)\setminus\{\ep\},$
$$\ep+\eta\not\in \bbbz(2\sgn(\ep;\a)\ep+k\d)+\sspan_\bbbz B=\sspan_\bbbz \Pi$$ which is a contradiction. Therefore,  $\ep\in\supp(B)$ and so there is $1\leq i\leq k$ with $\supp(B_i)\cap \supp(\a)\neq \emptyset.$

Now to the  contrary, assume there is  $1\leq j\leq k$ with $i\neq j$  such that  $\supp(B_i)\cap\supp(B_j)\cap  \supp(\a)\neq \emptyset.$  So there are $\b\in B_i$ and $\gamma\in B_j$ such that $\ep\in \supp(\b)\cap \supp(\gamma).$ Since  $(\b,\gamma)_*=0,$  Fact 1 implies that $\supp(\b)=\supp(\gamma)$ and  $|\supp(\b)|=|\supp(\gamma)|=2.$ Suppose $\supp(\b)=\supp(\gamma)=\{\ep,\eta\}.$ Since $\b-\a,\gamma-\a\not\in R,$
Facts~1,2 imply that
\begin{equation}
\label{fact2}
\sgn(\ep;\b)=-\sgn(\ep;\a)=\sgn(\ep;\gamma)\andd  \sgn(\eta;\b)=-\sgn(\eta;\a).
\end{equation}
Also using Lemma \ref{1 in 2}(b)(i),(b)(iii), we have
\begin{equation}\label{final}
\supp(\theta)\cap \{\ep,\eta\}=\emptyset\quad\quad(\theta\in\Pi\setminus\{\a,\b,\gamma\}).
\end{equation}
  We recall that $(n,m)\neq (1,1)$ and  choose  $\zeta\in \supp(\Pi)\setminus\{ \ep,\eta\},$  then $$\zeta+\sgn(\ep;\a)\ep=r\a+s\b+t\gamma+\sum_{\theta\in\Pi\setminus\{\a,\b,\gamma\}}r_\theta\theta$$ for some integers $r,s,t,r_\theta$ $(\theta\in \Pi\setminus\{\a,\b,\gamma\}).$ This together with (\ref{fact2}) and (\ref{final}) implies that $1=2r-s-t$ and $s=t$ which is a contradiction.

\medskip
\begin{centerline}
{\bf Using Step 4, without loss of generality, we assume $\boldsymbol{(\a,B_1)_*\neq \{0\}.}$}
\end{centerline}

\smallskip

\noindent{\bf Step 5.} ${B_1}$ is a base of ${S_1}:$
We must prove that $S_1\sub \sspan_{\bbbz^{\geq0}}B_1\cup \sspan_{\bbbz^{\leq0}}B_1.$ Suppose $\b\in S_1.$ We know that   $R\sub \sspan_{\bbbz^{\geq 0}}\Pi\cup \sspan_{\bbbz^{\leq 0}}\Pi.$ If\footnote{We remind that quasi-reflections do not necessarily preserves $R.$} $\b\in R,$ since $\b\in \sspan_\bbbz {B_1},$ we get that $\b\in \sspan_{\bbbz^{\geq 0}}{B_1}\cup \sspan_{\bbbz^{\leq 0}}{B_1}$ and so we are done.
Next assume $\b\in S_1\setminus R.$ Then  $\b\in \{\pm2\ep_i,\pm2\d_p\mid 1\leq i\leq m,1\leq p\leq n\}+\bbbz\d.$ The fact that  $(\a,B_1)\neq \{0\}$ together with Lemma \ref{unique} implies that  ${S_1}$ is not of simply-laced type. So  there is a long root $\gamma\in {B_1}.$ Assume $$\supp(\a)=\{\ep\},\;\;\supp(\gamma)=\{\zeta\}\andd \supp(\b)=\{\xi\}.$$
Since $(\a,B_1)_*\neq \{0\},$ by Lemma \ref{unique}, there is $\theta\in B_1$ with $|\supp(\theta)|=2$ and $\ep\in \supp(\theta).$ Also as $B_1$ is irreducible, again using Lemma \ref{unique}, there is $\theta'\in B_1$ with $|\supp(\theta')|=2$ and $\zeta\in \supp(\theta).$ Moreover, if
$\xi\neq \zeta,\ep,$ since $\b\in {S_1}=W_{B_1}{B_1},$ there is $\theta''\in B_1\setminus\{\a,\gamma\}$ with
 $|\supp(\theta'')|=2$ and $\xi\in \supp(\theta'').$
  These altogether imply that there is $t\in\{\pm 1\}$ and  $\rho\in(\{\ep_i,\d_p\mid 1\leq i\leq m, 1\leq p\leq n\}\setminus\{\xi\})+\bbbz\d$ with $t\xi+\rho\in {B_1}.$
Suuppose
\[\Pi=\{\a_1,\ldots,\a_{\ell+1}\}\quad\hbox{with}\quad\a_1=\a,\;\a_2=\gamma,\;\a_3:=t\xi+\rho\]  and
\[B_1=\{\a_2,\ldots,\a_{\ell'}\}\]
Without loss of generality we assume \[t=\sgn(\xi;\b).\]
Again using  Lemma \ref{unique}, we may assume $\ep,\xi\in\{\ep_i\mid 1\leq i\leq m\}$ and $\zeta\in\{\d_p\mid 1\leq p\leq n\}$ or $\ep,\xi\in \{\d_p\mid 1\leq p\leq n\}$ and $\zeta\in\{\ep_i\mid 1\leq i\leq m\}.$
So since $\b\not\in R,$ there are $r\in\{\pm1\},$ $k,k'\in\bbbz$ and $s\in\{0,1\}$ with
\[\b=2t\xi+2k\d+s\d\andd \gamma=2r\zeta+2k'\d+s\d.\]
Since $t\xi-\rho+2k\d+s\d,t\xi-r\zeta+(k-k')\d\in R,$ there are $r_i,t_i\in\bbbz$ ($1\leq i\leq \ell+1$) where  $t_i$'s (resp. $r_i$'s) are all non-positive or all non-negative with
\[t\xi-r\zeta+k\d-k'\d=\sum_{i=1}^{\ell+1}t_i\a_i\andd t\xi-\rho+2k\d+s\d= \sum_{i=1}^{\ell+1}r_i\a_i\]
and we have

\begin{align*}
  (1+r_3)\a_3+\sum_{3\neq i=1}^{\ell+1}r_i\a_i=&(t\xi-\rho+2k\d+s\d)+(t\xi+\rho)\\
  =&\underbrace{2t\xi+2k\d+s\d}_\b\\
  =&2(t\xi-r\zeta+k\d-k'\d)+2r\zeta+2k'\d+s\d\\
  =&(1+2t_2)\a_2+\sum_{2\neq i=1}^{\ell+1}2t_i\a_i.
\end{align*}

Since $\b\in S_1\sub \sspan_\bbbz B_1,$ we get
\begin{align*}
(1+r_3)\a_3+\sum_{3\neq i=2}^{\ell'}r_i\a_i=\b=(1+2t_2)\a_2+\sum_{i=3}^{\ell'}2t_i\a_i.
\end{align*}

If $\b\not\in \sspan_{\bbbz^{\geq0}} B_1\cup\sspan_{\bbbz^{\leq0}} B_1,$ using the same argument as in the proof of  Proposition \ref{with-long},  we get  $r_3=0$ and so $1=2t_3,$ a contradiction.

\noindent{\bf  Step 6.} $k=1:$ Suppose $k\geq 2.$ Recall that $S_1$ is an irreducible finite root system with base $B_1$ and that   $\{0\}\neq (\a_1,B_1)_*\sub (\a_1,S_1)_*.$ Using Lemma \ref{types}, we have the following cases:

\noindent \underline{\it Case 1.} $S_1$ is of type $D_{\frak{r}}$ $(\frak{r}\geq 4):$
By Corollary \ref{cor1}(i) and Step 5,  there are   $$r\in\{0,1\},\;\;k^*\in\bbbz\andd\zeta_1,\ldots,\zeta_{ \frak{r}}\in \{\pm\ep_{\frak{r}},\pm\d_p\mid 1\leq t\leq n,1\leq p\leq m\}+\bbbz\d$$  with   $\supp(\zeta_i)=\{\dot\zeta_i\}$  such that  $\supp(\a)=\{\dot \zeta_1\},$  $\dot\zeta_i$'s are distinct,
    $$\a=-2\zeta_1+r\d\andd B_1=\{\gamma_i:=\zeta_i-\zeta_{i+1},\gamma_{\frak{r}}:=\zeta_{\frak{r}-1}+\zeta_{\frak{r}}+k^*\d\mid 1\leq i\leq \frak{r}-1\}.$$
 Now we have  $$(r+k^*)\d=\gamma+2\a_1+\cdots+2\gamma_{\frak{r}-2}+\a_{\frak{r}-1}+\a_{\frak{r}}.$$ Since $\{\a,\gamma_1,\ldots,\gamma_{\frak{r}}\}$ is linearly independent, $k^*+r\neq 0,$ so $\d\in \sspan_\bbbz B_1\cup \{\a\}.$ Next fix a root  $\rho\in B_2.$ If $\d\in \sspan_{\bbbz^{\geq0}} B_1\cup \{\a\},$ then $-\rho+2\d $ is an element of $R\cap (\sspan_{\bbbz^{<0}} B_2+\sspan_{\bbbz^{\geq0}} B_1\cup \{\a\})$
which is a contradiction. Similarly,    if $\d\in \sspan_{\bbbz^{\leq0}} B_1\cup \{\a\},$ we have $\rho+2\d \in R\cap (\sspan_{\bbbz^{>0}} B_2+\sspan_{\bbbz^{\leq0}} B_1\cup \{\a\})$
which is again a  contradiction.

\smallskip

\noindent \underline{\it Case 2.} $S_1$ is of type $C_{\frak{r}}$ $(\frak{r}\geq 2).$ We get a contradiction using the same argument as in Case 1.

\smallskip

\noindent\underline{\it Case 3.} $S_1$ is of type $A_{\frak{r}}~(\frak{r}\neq3).$
By Corollary \ref{cor1}(iii) and Step 5, there are   $$r\in\{0,1\},\;\;\zeta_1,\ldots,\zeta_{ \frak{r}+1}\in \{\pm\ep_t,\pm\d_p\mid 1\leq t\leq n,1\leq p\leq m\}+\bbbz\d$$  with   $\supp(\zeta_i)=\{\dot\zeta_i\}$  such that  $\supp(\a)=\{\dot \zeta_1\},$  $\dot\zeta_i$'s are distinct,
    $$\a=-2\zeta_1+r\d\andd B_1=\{\a_i:=\zeta_i-\zeta_{i+1}\mid 1\leq i\leq \frak{r}\}.$$
 If $2\leq j\leq \frak{r}+1$ and  $\dot\zeta_{j}\in \supp(B_i)$ for some $i\neq 1,$ then there  is $\b\in B_i$ with $\dot\zeta_j\in \supp(\b).$ Since $(B_1,B_i)_*=\{0\},$  we have $(\a_{j-1},\b)_*=(\a_{j},\b)_*=0$ and so by Fact 1, we have $\supp(\a_{j-1})=\supp(\b)=\supp(\a_j)$ which is a contradiction. Therefore, contemplating Step 4. we have
 \begin{equation}\label{cap}
 \supp(B_1)\cap \cup_{i=2}^k\supp(B_i)=\emptyset.
 \end{equation}
But there are $r_1,\ldots,r_{\frak{r}+1}\in\bbbz$ with  \[\d=r_1(-2\zeta_1+r\d)+\sum_{i=2}^{\frak{r}+1}r_i(\zeta_{i-1}-\zeta_i)+\sum_{\a\in \cup_{i\geq 2} B_i}k_\a\a\] So  using  (\ref{cap}), $2r_1=r_2=\cdots=r_{\frak{r}+1}=0,$ in other words, $\d\in \sspan_{\bbbz}\cup_{i=2}^kB_i$ which is a contradiction as the form $\fm_*$ is nondegenerate on $\sspan_\bbbr S.$

\smallskip

\noindent\underline{\it Case 4.} $S_1$ is  of type $A_3:$  Using Corollary \ref{cor1}(iv) together with the same argument as in Cases 1 and 3, we get a contradiction. This completes the proof.

\noindent{\bf Step 7.} If $\ell\neq 3,$ then  $S$ is not of type $A:$ Since $S=S_1$ is of rank $\ell=m+n,$ if it is of type $A,$ by Corollary \ref{cor1}(iii) and Step 5,  there are $\ell+1$ elements $\zeta_1,\ldots,\zeta_{\ell+1}$ of $\{\pm\ep_i,\pm\d_p\mid 1\leq i\leq m,1\leq p\leq n\}$ with $\zeta_i\neq \pm\zeta_j$ if $i\neq j$ which is a  contradiction.

\medskip

Steps 3,5,6 and 7 together with Lemma \ref{types} complete the proof.
\qed

\medskip

\begin{lem}\label{b1}
 Suppose that $B$ is  a base of $R$ without long roots, then $B$ is of the form $\pm B^1.$
 Set
 \[
 \xi_i:=\left\{\begin{array}{ll}
\d_i& 1\leq i\leq n,\\
\ep_{i-n}& n+1\leq i\leq \ell
\end{array}\right.
\] and
\begin{align*}
\Pi:=\{\d-(\xi_1+\xi_2),\xi_i-\xi_{i+1},\xi_{\ell-1}+\xi_\ell\mid 1\leq i\leq \ell-1\}.
\end{align*}
Then $\Pi$ is a base of $R$ and   each base of $R$ which does not contain  long roots is conjugate  with $\pm\Pi$ under quasi-Weyl group.
\end{lem}
\pf
That $\Pi$ is a base of $R$ follows from Proposition \ref{bases-last}. Suppose that $B$ is  a base of $R$ which does not contain long roots.
 Set
\begin{align*}
S:=&\pm\{0,\ep_i\pm\ep_j,\d_p\pm\d_q,\ep_i\pm\d_p\mid 1\leq i\neq j\leq m,1\leq p\neq q\leq n\}+\bbbz\d\\
=&\{0,\pm\xi_i\pm\xi_j\mid 1\leq i\neq j\leq \ell\}+\bbbz\d.
\end{align*}
Then $(\sspan_\bbbr S,\fm_*,S)$ is the root system of $D_\ell^{(1)}$ and $\Pi$ is a base of $R.$ Since $B$ is a subset of $S,$  it is  a base of $S.$ So  there is an element $\omega$ of the Weyl group of $S$ such that $\omega(\Pi)= r B$ for some $r\in\{\pm1\}.$ Setting
\[\theta_i:=r\omega(\xi_i)\in\{\pm \xi_j\mid 1\leq j\leq \ell\}+\bbbz\d\quad\quad(1\leq i\leq \ell),\]
we have
$$B=\{r\d-(\theta_1+\theta_2),\theta_i-\theta_{i+1},\theta_{\ell-1}+\theta_\ell\}.$$
Since  $$2\theta_1-r\d=(\theta_1-\theta_2)-(r\d-(\theta_1+\theta_2)),2\theta_\ell=(\theta_{\ell-1}+\theta_\ell)-(\theta_{\ell-1}-\theta_{\ell}) \in B-B,$$ we have  $$2\theta_1-r\d, 2\theta_\ell\not\in R$$ and so $$\supp(\theta_1)\sub\{\d_p\mid 1\leq p\leq n\}\andd\supp(\theta_\ell)\sub\{\ep_i\mid 1\leq i\leq m\}. $$

 This completes the proof as the Weyl group of $S$ is a subgroup of the quasi-Weyl group of $R.$ \qed

\begin{Pro}\label{bases}
\begin{itemize}
\item[(i)]   If $\Pi$ is  a base of $R$ without long roots, then $\Pi$ is of the form $\pm B^1.$
\item[(ii)]
If $\Pi$ is a base of  $R$  having a unique  long root which belongs to  $\{\pm2\d_p\mid 1\leq p\leq n\}+2\bbbz\d,$ then $\Pi$ is of the form $\pm B^2.$
\item[(iii)] If $\Pi$ is a base of $R$   having a unique  long root which belongs to  $\{\pm2\ep_i\mid 1\leq i\leq m\}+2\bbbz\d+\d,$ then it is of the form $\pm B^3.$
\item[(iv)]  If $\Pi$ is a base of $R$ having two  long roots, then it  is of the form $\pm B^4.$
\end{itemize}
\end{Pro}
\pf
(i) See Lemma \ref{b1}.

(ii) Using   Theorem \ref{finite} and  Corollary \ref{cor1}(i),(iv), we find
$k^*\in\bbbz$ and $\theta_1,\ldots,\theta_{ \ell}\in \{\pm\ep_t,\pm\d_p\mid 1\leq t\leq m,1\leq p\leq n\}+\bbbz\d$  with   $\supp(\theta_i)=\{\upsilon_i\}$  such that   $\upsilon_i$'s are distinct and
    $$\Pi=\{-2\theta_1,\theta_i-\theta_{i+1},\theta_{\ell-1}+\theta_{\ell}+k^*\d\mid 1\leq i\leq \ell-1\}.$$
  Since $\d\in \sspan_\bbbz\Pi,$  it follows that $k^*=\pm 1.$ Also as
  \[2\theta_{\ell}+k^*\d=(\theta_{\ell-1}+\theta_{\ell}+k^*\d)-(\theta_{\ell-1}-\theta_{\ell})\in\Pi-\Pi,\] we get that  $2\theta_{\ell}+k^*\d\not\in R$ and so $\upsilon_\ell\in \{\d_p\mid 1\leq p\leq n\}.$

(iii) Using   Theorem \ref{finite} and  Corollary \ref{cor1}(i),(iv), we find
$k^*\in\bbbz$ and $\theta_1,\ldots,\theta_{ \frak{r}}\in \{\pm\ep_t,\pm\d_p\mid 1\leq t\leq n,1\leq p\leq m\}+\bbbz\d$  with   $\supp(\theta_i)=\{\upsilon_i\}$  such that   $\upsilon_i$'s are distinct and
    $$\Pi=\{-2\theta_1+\d,\theta_i-\theta_{i+1},\theta_{\ell-1}+\theta_{\ell}+k^*\d\mid 1\leq i\leq \ell-1\}.$$
  Since $\d\in \sspan_\bbbz\Pi,$  it follows that $k^*\in\{0,-2\}.$
  Also as
  \[2\theta_{\ell}+k^*\d=(\theta_{\ell-1}+\theta_{\ell}+k^*\d)-(\theta_{\ell-1}-\theta_{\ell})\in\Pi-\Pi,\] we get that  $2\theta_{\ell}+k^*\d\not\in R$ and so $\upsilon_\ell\in \{\ep_i\mid 1\leq i\leq m\}.$
  Therefore, if $k^*=0,$ we are done. But if $k^*=-2,$ setting
  $\theta'_i:=\theta_i-\d$ $(1\leq i\leq \ell),$  we have
  $$\Pi=\{-2\theta'_1-\d,\theta'_i-\theta'_{i+1},\theta'_{\ell-1}+\theta'_{\ell}\mid 1\leq i\leq \ell-1\}$$ which is of the form $- B^3.$

 (iv) Using   Theorem \ref{finite} and  Corollary \ref{cor1}(ii), we find
$k^*\in\bbbz$ and $\theta_1,\ldots,\theta_{ \frak{r}}\in \{\pm\ep_t,\pm\d_p\mid 1\leq t\leq n,1\leq p\leq m\}+\bbbz\d$  with   $\supp(\theta_i)=\{\upsilon_i\}$  such that   $\upsilon_i$'s are distinct and
    $$\Pi=\{-2\theta_1,\theta_i-\theta_{i+1},2\theta_{\ell}+k^*\d\mid 1\leq i\leq \ell-1\}.$$ Since $\d\in \sspan_\bbbz\Pi,$  it follows that $k^*\in\{\pm 1\}.$ This completes the proof.\qed

\begin{Thm}
Suppose $R=A(2m-1,2n-1)^{(2)}$ for some positive integers $m,n$ with $(m,n)\neq (1,1).$ Then we have the following:
\begin{itemize}
\item[(i)] Each base of $R$ is of the form $\pm B^j$ for some $1\leq j\leq 4.$
\item[(ii)]
Set \[
 \xi_i:=\left\{\begin{array}{ll}
\d_i& 1\leq i\leq n,\\
\ep_{i-n}& n+1\leq i\leq \ell
\end{array}\right.
\]and
\begin{align*}
\Pi^1:=&\{\d-(\xi_1+\xi_2),\xi_i-\xi_{i+1},\xi_{\ell-1}+\xi_\ell\mid 1\leq i\leq \ell-1\},\\
\Pi^2:=&\{-2\xi_2,\xi_{i}-\xi_{i+1},\xi_\ell-\xi_1,\xi_{\ell}+\xi_1+\d\mid 2\leq i\leq \ell-1\}\quad(\hbox{if $n\geq 2$}),\\
\Pi^3:=&\{-2\xi_\ell+\d,\xi_\ell-\xi_1,\xi_i-\xi_{i+1},\xi_{\ell-2}+\xi_{\ell-1}\mid 1\leq i\leq \ell-2\}\quad(\hbox{if $m\geq 2$}),\\
\Pi^4:=&\{-2\xi_1,\xi_i-\xi_{i+1},2\xi_\ell+\d\mid 1\leq i\leq \ell-1\}.
\end{align*} Then $\Pi^1,\ldots,\Pi^4$ are bases of $R$ and  if  $\Pi$ is a base of $R,$ then there is $w\in W$ and $1\leq i\leq 4$ such that $\Pi=\pm w(\Pi^i).$
\item[(iii)] Each base of $R$ is an S-base.
\end{itemize}
\end{Thm}
\pf
(i) See Proposition \ref{bases}.

(ii) follows from part (i), Lemma \ref{b1} and   Proposition \ref{bases-last} together with Corollary \ref{cor2}.

(iii) Set
\[\dot R:=\{0,\pm2\ep_i,\pm2\d_p,\pm\ep_i\pm\ep_j,\pm\d_p\pm\d_q,\pm\ep_i\pm\d_p\mid 1\leq i\neq j\leq m,1\leq p\neq q\leq n \}\] and assume
$\Pi$ is a base of $R.$ Using  part (i) together with Lemma  \ref{b1} and Proposition   \ref{bases-last}, for each $\a\in \dot R,$ there is a positive integer $N$ such that either  \[R\cap (\a+\bbbz^{\geq N})\sub R^+(\Pi) \andd R\cap (\a+\bbbz^{\leq N})\sub R^-(\Pi)\] or
\[R\cap (\a+\bbbz^{\leq N})\sub R^+(\Pi) \andd R\cap (\a+\bbbz^{\geq N})\sub R^-(\Pi).\]
 This implies that if $\Pi$ and $\Pi'$ are two bases of $R,$  then either $R^+(\Pi)\cap R^{-}(\Pi')$ or $R^+(\Pi)\cap R^{+}(\Pi')$ is a finite set. In particular, if $\Pi$ is the standard base of $R,$  for each base $\Pi'$ of $R,$ either
$R^+(\Pi)\cap R^-(\Pi')$ or $R^+(\Pi)\cap R^{+}(\Pi')$ is a finite set.  Now the result follows using Lemma \ref{new} as the standard base is in fact an S-base.\qed

\end{document}